\documentclass{amsart}
\usepackage{amsmath,amsthm,amssymb}
\usepackage{lipsum} 
\usepackage[utf8]{inputenc}
\usepackage{esint}
\usepackage{amsfonts} 
\usepackage{xcolor}
\usepackage{hyperref}
\usepackage{bbold}
\usepackage{enumerate}
\usepackage{pgfplots}
\usepackage{siunitx}
\pgfplotsset{width=11cm,compat=1.9}
\usepgfplotslibrary{fillbetween}
\usetikzlibrary{patterns}
\makeatletter
\numberwithin{equation}{section}
\theoremstyle{plain}
        \newtheorem{theorem}{Theorem}[section]
        
        \newtheorem{lemma}[theorem]{Lemma}

        \newtheorem{definition}[theorem]{Definition} 
        \newtheorem{remark}[theorem]{Remark}  
         
        \newtheorem*{claim*}{Claim}  
         
\newtheorem*{theorem*}{Theorem}
\newtheorem*{definition*}{Definition}
\newtheorem*{proposition*}{Proposition}
\usepackage{graphicx}
\let\oldmarginpar\marginpar
\renewcommand\marginpar[1]{\-\oldmarginpar[\raggedleft\footnotesize #1]
{\raggedright\footnotesize #1}}

\def\dx{\,{\rm d}x\,}

\def\dt'{\,{\rm d}t'\,}

\renewcommand\div{\text{div\,}}
\newcommand \loc {\text{loc}}

\newcommand{\R}{\mathbb{R}}

\newcommand{\PP}{\mathbb{P}}
\newcommand{\N}{\mathbb{N}}

\newcommand{\cD}{\mathcal{D}}

\newcommand{\cF}{\mathcal{F}}
\newcommand{\cS}{\mathcal{S}} 
\newcommand\Lip{\mathrm{Lip}}

\renewcommand{\d}{\partial} 
\renewcommand\dt{\mathrm{dt}}
\DeclareMathOperator{\di}{div}

\title[Two dimensional Boussinesq equations in general Sobolev spaces]{On the two-dimensional Boussinesq equations with temperature-dependent thermal and viscosity diffusions in general Sobolev spaces}
\author[Z. He]{Zihui He} 
\author[X. Liao]{Xian Liao} 
\address[Z. He and X. Liao]
{Institute of Analysis, Karlsruhe Institute of Technology\\
Englerstraße 2, 76131 Karlsruhe, Germany.}
\email{zihui.he@kit.edu, xian.liao@kit.edu}
\begin{document}
\subjclass[2020]{35Q35, 76D03} 
\keywords{Boussinesq equations, variable diffusion coefficients, existence, uniqueness, regularity} 
\maketitle
\begin{abstract}
We study the existence, uniqueness as well as regularity issues for the two-dimensional incompressible Boussinesq equations with temperature-dependent thermal and viscosity diffusion coefficients in general Sobolev spaces. 
  The optimal regularity exponent ranges are considered.
\end{abstract}
\bigskip 
\section{Introduction}
In the present paper we consider the two-dimensional incompressible Boussinesq equations
\begin{equation}\label{BS}
\left\{
\begin{aligned}
&\d_t\theta+ u\cdot \nabla_x \theta-\di_x(\kappa \nabla_x \theta)=0,\\
&\d_t u+u\cdot \nabla_x u-\di_x(\mu  S_x u)+\nabla_x\Pi=\beta \theta \vec{e_2},\\
&\di_x u=0, 
\end{aligned}
\right.
\end{equation}
where $(t,x)\in [0,\infty)\times \R^2$ denote the time and space variables respectively.
The unknown temperature function $\theta=\theta(t,x):[0,\infty)\times\R^2\to \R$ satisfies the parabolic-type equation $\eqref{BS}_1$, and the unknown velocity vector field $u=u(t,x): [0,\infty)\times\R^2\to\R^2$ together with the unknown pressure term $\Pi=\Pi(t,x):[0,\infty)\times\R^2\to \R$ satisfies  the incompressible Navier-Stokes type equations $\eqref{BS}_2-\eqref{BS}_3$ respectively.
We are going to study   the well-posedness and regularity   problems for  the Boussinesq system \eqref{BS} together with the initial data 
\begin{equation}\label{BS-IC}
(\theta,u)\mid_{t=0}=(\theta_0,u_0).
\end{equation}

We write $x=\begin{pmatrix}x_1\\ x_2\end{pmatrix}\in \R^2$ with $x_1, x_2$ denoting the horizontal and vertical components respectively.
Let $u=\begin{pmatrix}u^1\\ u^2\end{pmatrix}: [0,\infty)\times\R^2\to\R^2$,
and  let
$$\frac12 S_xu:=\frac12 (\nabla_x u+(\nabla_x u)^T),\hbox{ with }\nabla_x u=(\partial_{x_j}u^i)_{1\leq i,j\leq 2}$$
denote the symmetric deformation  tensor in the second equation $\eqref{BS}_2$ above.
The vector field $\vec{e_2}$ denotes the unit vector in the vertical direction: $\vec{e_2}=\begin{pmatrix}0\\1\end{pmatrix}$,
and $\beta\theta\vec{e_2}$ stays for the buoyancy force, with the constant parameter $\beta>0$ denoting the  thermodynamic dilatation coefficient which will be assumed to be $1$ in the following context for simplicity.

We consider the cases when the heat diffusion and the viscosity in the fluids are sensitive to the change of temperatures, that is,  the thermal diffusivity $\kappa$  and the viscosity coefficient $\mu$ may depend on the temperature function $\theta$ smoothly as follows
\begin{equation}\label{ab}
\kappa=a(\theta),\quad \mu=b(\theta), \quad \hbox{ with }a\in C^1_b (\R; [\kappa_\ast, \kappa^\ast]),\, b\in C^1_b (\R; [\mu_\ast, \mu^\ast]),
\end{equation}
where $\kappa_\ast\leq \kappa^\ast$, $\mu_\ast\leq \mu^\ast$ are positive constants.
We will not assume any smallness conditions on   $\kappa^\ast-\kappa_\ast$ or $\mu^\ast-\mu_\ast$, and large variations in these diffusivity coefficients are permitted.  

\medskip

The Boussinesq system \eqref{BS} arises from the zero order approximation to the corresponding inhomogeneous hydrodynamic systems, which are nonlinear coupling between the Navier-Stokes equations or Euler equations and the thermodynamic equations for the temperature or density functions: The Boussinesq approximation \cite{Boussinesq} ignores density differences expect when they appear in the buoyancy term.
They are common geophysical models describing the dynamics from  large scale atmosphere and ocean flows to solar and plasma inner convection, where density stratification is a typical feature \cite{Gill, Majda}.

The temperature differences in the inhomogeneous fluids may cause \emph{density gradients}.
When the thermodynamical coefficients such as the heat-conducting coefficients and the viscosity coefficients are assumed to be constant in the Boussinesq approximation (i.e. $\kappa$, $\mu$ are constants in \eqref{BS}),  density gradients  influence the motion of the flows only through the buoyancy force, which may lead to finite time singularity in the flows (the formation of the finite time singularity is sensitive to the thermal and viscous dissipation and see Subsection \ref{subs:ref} below for more references on this topic).

However, the temperature variations do influence the thermal conductivity and the viscosity coefficients effectively, even for simple fluids such as pure water \cite[Section 6]{Lide}\footnote{The absolute viscosity of the water under  nominal atmospheric pressure in units of millipascal seconds is given by $1.793$ (\SI{0}{\celsius}), $0.547$ (\SI{50}{\celsius}),  $0.282$ (\SI{100}{\celsius}) respectively \cite[Page 6-186]{Lide}.
The thermal conductivity of the water under nominal atmospheric pressure in units of  watt per meter kelvin is given by  $0.5562$ (\SI{0}{\celsius}), $0.6423$ (\SI{50}{\celsius}),  $0.6729$ (\SI{100}{\celsius}) respectively \cite[Page 6-214]{Lide}.}
\footnote{It is common to adapt the exponential viscosity law $\mu(T)=C_1\exp(C_2/(C_3+T))$ and quasi-constant heat conductivity law $\kappa(T)=C_4$ for the liquids, while the viscosity law $\mu(T)=(\mu(T_m)) \frac{T}{T_m}\frac{T_m+C_5}{T+C_6}$ and the thermal conductivity law $\kappa(T)=C_6\mu(T)$  for the gases, where $T$ denotes the absolute temperature, $T_m$ denotes the reference temperature and $C_j$, $1\leq j\leq 6$ are positive constants \cite[I]{Perez}.  }.
In many applications in the engineering one also aims for  effective thermal conductivities in building thermal energy storage materials \cite{Gaedtke}.
Therefore in many physical models density gradients would influence the motion of the fluids not only through buoyancy force, but also through the variations of the diffusion coefficients.
It is then interesting to study the wellposedness and regularity problems of the Boussinesq system \eqref{BS}-\eqref{ab}.
 
\subsection{Known results}\label{subs:ref}
The wellposedness and regularity problems on the   two  - dimensional Boussinesq equations have attracted considerable attention from the PDE community. 
Many interesting mathematical results have been established in the past two decades, mainly in the cases with constant thermal diffusivity coefficient $\kappa$ and viscosity coefficient $\mu$:
\begin{equation}\label{BS:constant}
\left\{
\begin{aligned}
&\d_t\theta+ u\cdot \nabla_x \theta-\kappa \Delta_x \theta=0,\\
&\d_t u+u\cdot \nabla_x u-\mu \Delta_x u+\nabla_x\Pi= \theta \vec{e_2},\\
&\di_x u=0,\\
&(\theta,u)\mid_{t=0}=(\theta_0,u_0).
\end{aligned}
\right.
\end{equation}

If $\kappa=\mu=0$, the two-dimensional inviscid Boussinesq equations \eqref{BS:constant} can be compared with the three-dimensional incompressible axisymmetric Euler equations with swirl, where the buoyancy force corresponds to the vortex stretching mechanism. 
The local-in-time wellposedness   as well as some blowup criteria have been well known for decades, see e.g.   \cite{ChaeNam, Danchin, EShu}.
We mention that an (improved) lower bound for the lifespan which tends to infinity as the initial temperature 
tends to a constant (and correspondingly, as the initial swirl
 tends to zero for the $3D$ axisymmetric Euler equations) was given in \cite{Danchin}. 
The outstanding global-in-time regularity problem has been solved recently in \cite{Elgindi:Jeong}, where  examples of finite-energy strong solutions which become singular in finite time have been given  (see examples of  finite time singularity solutions for $3D$ axisymmetric Euler equations in  \cite{Elgindi}).

If $\kappa>0$, $\mu>0$, on the contrary,  the convection terms can be controlled thanks to the strong diffusion effects, and the global-in-time existence and regularity results can be established (see e.g.   \cite{Cannon:DiBenedetto}).
Particular interests then raised if only partial dissipation is present, that is, either $\kappa=0$ whereas $\mu>0$ or $\kappa>0$ whereas $\mu=0$ (see e.g. H.K. Moffatt's list of the 21st Century PDE problems \cite{Moffatt}).
The global-in-time results continue to hold, thanks to  a priori estimates in the $L^p$-framework as well as the sharp Sobolev embedding inequality in dimension two with a  logarithm correction,  which help the partial diffusion terms to control the demanding term $\partial_{x_1}\theta$ successfully (see \cite{Chae, Hou:Li}   and see  \cite{Hmidi:Keraani} for less regular cases). 
Further developments were made for horizontal dissipation cases (see e.g. \cite{Danchin:Paicu}), for vertical dissipation cases (see e.g. \cite{Cao:Wu}), and for the fractional dissipation cases (see e.g. \cite{Hmidi:Keraani:Rousset1, Hmidi:Keraani:Rousset2}).
See the review notes \cite{Wu} and the references therein for more interesting results and sketchy proofs.  

There also have been  remarkable progresses in solving   the two dimensional Boussinesq equations \eqref{BS}-\eqref{ab} when  the thermal and viscosity diffusion coefficients $\kappa, \mu$ are variable and depend smoothly on the unknown temperature function $\theta$.
In the variational formulation framework, the global-in-time existence of a solution of \eqref{BS}-\eqref{ab} has been established in \cite{Duvaut:Lions} (see \cite{Feireisl:Malek} for a similar formulation of \eqref{BS}-\eqref{ab}) for the motion of the so-called Bingham fluid (as a non-Newtonian fluid), where $\kappa$ is a positive constant, $\beta=0$ and $\mu$ depends not only on $\theta$ but also on  $Su/|Su|$. 
The Boussinesq-Stefan model has been investigated in \cite{Rodrigues}, where the phase transition was taken into account.
The global-in-time existence as well as the uniqueness of the solutions for  \eqref{BS}-\eqref{ab} have been shown in \cite{DiazGaliano, Goncharova, Lorca:Boldrini} under Dirichlet boundary conditions and in \cite{Perez} under generalized outflow boundary conditions. We remark that the resolution of the nonhomogeneous Boussinesq system under more physical boundary conditions (e.g. with Dirichlet boundary conditions only on the inflow part of the boundary while with no prescribed assumptions on the outflow part) remains unsolved.

S. Lorca and J. Boldrini \cite{Lorca:Boldrini} (see also \cite{DiazGaliano, Goncharova}) studied the initial-boundary value problem of the Boussinesq system \eqref{BS}-\eqref{ab} in dimension two and three  under the initial condition \eqref{BS-IC} and Dirichlet boundary conditions,
and obtained a global-in-time weak solution 
$$(\theta,u)\in(L^\infty_\loc([0,\infty); L^2(\Omega)) )^3$$
as well as a local-in-time unique strong solution  
\begin{equation}\label{LB}
(\theta,u)\in L^\infty_\loc([0,\infty); H^2(\Omega))\times (L^\infty_\loc([0,\infty); H^{1}(\Omega))^2.
\end{equation}
The outstanding global-in-time existence and uniqueness results of the  smooth solutions 
\begin{equation}\label{WZ}(\theta,u)\in (L^\infty_\loc([0,\infty); H^s(\R^2))\cap L^2_\loc([0,\infty); H^{s+1}(\R^2)))^3,\quad s>2\end{equation}
 have been successfully established by C. Wang and Z. Zhang \cite{Wang:Zhang}, which affirms the propagation of high regularities (without finite time singularity) of the two dimensional Boussinesq flow in the presence of viscosity variations (see \cite{Sun:Zhang} for the case $s=2$).

We remark that it is still not clear whether there will be finite time singularity for the two dimensional Boussinesq flow \eqref{BS}-\eqref{ab} in the presence of viscosity variations while no heat diffusion (i.e. $\kappa=0$, $\mu=\mu(\theta)$), and we mention a recent work \cite{Abidi:Zhang} toward this direction in the case of less heat diffusion (with $\div(\kappa\nabla\theta)$ replaced by $(-\Delta)^{1/2}$) and the small viscosity variation assumption: $|\mu-1|\leq \varepsilon$.
A closely related question would pertain to the global-in-time wellposedness problem of the two-dimensional inhomogeneous incompressible Navier-Stokes equations with density-dependent viscosity coefficient 
\begin{equation}\label{NS}
\left\{
\begin{aligned}
&\d_t\rho+u\cdot\nabla_x\rho=0,\\
&\partial_{t}(\rho u)+\operatorname{div}_x(\rho u\otimes u)-\operatorname{div}_x(  \mu S_xu)+\nabla_x \Pi=0, \\
&\operatorname{div}_x u=0, \\
&(\rho, \, \rho u)\big|_{t=0}=(\rho_{0}, m_{0}).
\end{aligned}
\right.
\end{equation}
The global-in-time existence results of weak solutions of \eqref{NS} (see e.g. \cite{AKN, Lions}) as well as the local-in-time well-posedness results (see e.g. \cite{LaSo}) have been well known,  while the global-in-time regularities still remain open (see e.g. \cite{AbidiZhang, Desjardins} for some interesting results  under the  assumption on the weak inhomogeneity).

 To the best of our knowledge, there are no global-in-time regularity propagation results  by the two-dimensional Boussinesq flow with temperature-dependent diffusion coefficients  \eqref{BS}-\eqref{BS-IC}-\eqref{ab} in the low regularity regime 
 $$H^s,\quad s<2,$$
 or in the general Sobolev setting
 $$\theta_0\in H^{s_\theta}_x(\R^2),\quad u_0\in (H^{s_u}_x(\R^2))^2$$ with different regularity indices $s_\theta$ and $s_u$. 
 In this paper we are going to investigate the existence, uniqueness as well as the regularity problems in these general Sobolev functional settings.

To conclude this subsection let us just mention some recent interesting progresses on the stability of the stationary shear flow solutions (together with the corresponding striated temperature function) to the Boussinesq equations \eqref{BS:constant}, with full dissipation or partial dissipation, in e.g. \cite{TaoWuZhaoZheng, Zillinger} and references therein. 
It should also be interesting to investigate the stability of the   stationary striated solutions of the Boussinesq equations with variable diffusion coefficients \eqref{BS}.



\subsection{Main results}
We are going to show the global-in-time existence of  weak solutions to the Cauchy problem for the Boussinesq system \eqref{BS}-\eqref{BS-IC}-\eqref{ab} in the whole two-dimensional space $\R^2$  under the low-regularity initial condition
 $(\theta_0, u_0)\in L^2(\R^2)\times (L^2(\R^2))^2.$ 
The uniqueness result holds true if the initial temperature function
becomes smoother $(\theta_0, u_0)\in H^1(\R^2)\times (L^2(\R^2))^2.$
Finally we will establish the global-in-time regularity of the solutions  in the general Sobolev setting 
$ (\theta_0, u_0)\in H^{s_\theta}(\R^2)\times (H^{s_u}(\R^2))^2\subset H^1(\R^2)\times (L^2(\R^2))^2$ with the restriction $ s_u-1\leq s_\theta\leq s_u+2$.
These regularity exponent ranges are  optimal for the existence, uniqueness and regularity results respectively, by view of the formulations   of the Boussinesq equations \eqref{BS} with temperature-dependant diffusion coefficients (see Remark \ref{rem:optimality} below for more details).

\begin{figure}
\centering
\begin{tikzpicture}
\begin{axis}[
    axis lines = middle,
    xlabel = $s_\theta$,
    ylabel = {$s_u$},
    xmin=0, xmax=5,
    ymin=0, ymax=5,]
    
\addplot [
    name path=A,
    domain=2:5, 
    samples=2,]{x-2};

\addplot [
    name path=C,
    domain=1:2, 
    samples=2,]{0};

\addplot [
name path=B,
    domain=1:4, 
    samples=2, 
    ]
    {x+1};

    

 
 \addplot [gray]  fill between [of=A and B];
 
  \addplot [gray]  fill between [of=C and B];
 






\end{axis}
\end{tikzpicture}
\caption{Admissible regularity exponents} \label{fig:s}
\end{figure}
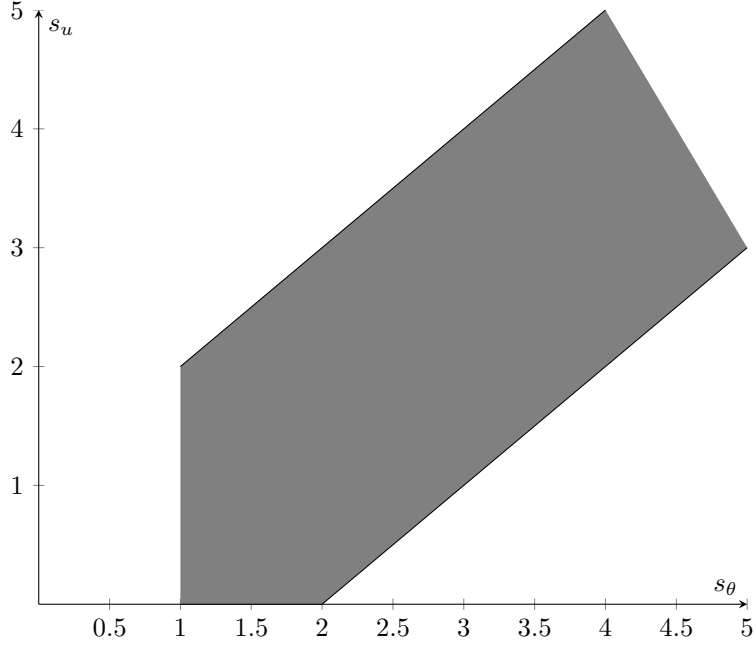

 We first define the weak solutions as follows.
\begin{definition}[Weak solutions]\label{def:weak} We say that a pair 
 $(\theta,u)$ is a weak solution of the Boussinesq equations \eqref{BS}-\eqref{ab} with the given initial data $(\theta_0,u_0)\in   (L^2(\R^2))^3$ if the following statements hold:
\begin{itemize}
    \item The temperature function 
    $$\theta=\theta(t,x)\in C([0,\infty); L^2_x(\R^2))\cap L^2_\loc([0,\infty); H^1_x(\R^2))$$
     satisfies the initial condition $\theta|_{t=0}=\theta_0$, the energy equality 
    \begin{equation}\label{EE:theta}
    \frac12\|\theta(T,\cdot)\|_{L^2_x(\R^2)}^2+\int^T_0\int_{\R^2} (\kappa|\nabla\theta|^2)(t, x)  \dx \dt 
    =\frac12\|\theta(0,\cdot)\|_{L^2_x(\R^2)}^2, 
    \end{equation}
   for all positive times $T>0$,  and the equation 
     \begin{equation}\label{BS:theta}
     \d_t\theta+u\cdot\nabla\theta-\div_x(\kappa\nabla\theta)=0
     \end{equation}
  in $L^2_\loc([0,\infty); H^{-1}_x(\R^2))$.
    \item The velocity vector field 
    $$u=u(t,x)\in C([0,\infty); (L^2_x(\R^2))^2)\cap L^2_\loc([0,\infty); (H^1_x(\R^2))^2)$$
     satisfies  the initial condition  $u|_{t=0}=u_0$, the divergence-free condition $\div_x u=0$, the energy equality
    \begin{equation}\label{EE:u}\begin{split}
    &\frac12\| u(T,\cdot)\|_{L^2_x(\R^2)}^2+\frac12\int^T_0\int_{\R^2} (\mu|S u|^2)(t, x)\dx \dt 
    \\
&    =\frac12\|u(0,\cdot)\|_{L^2_x(\R^2)}^2+\int^T_0\int_{\R^2} (\theta u_2)(t,x) \dx \dt ,\quad \forall T>0,
  \end{split}  \end{equation}
      and the equation 
      \begin{equation}\label{BS:u}
      \d_t u+u\cdot \nabla_x u-\di_x(\mu  S_x u)+\nabla_x\Pi= \theta \vec{e_2}
      \end{equation}
 in  $L^2_\loc([0,\infty); (H^{-1}_x(\R^2))^2 )$ for some scalar function $\Pi\in L^2_\loc([0,\infty)\times\R^2)$ with $\nabla\Pi\in L^2_\loc([0,\infty); (H^{-1}_x(\R^2))^2 ) $ and $\int_{B_1}\Pi\dx=0$ a.e. $t$ (with $B_1$ denoting the unit disk in $\R^2$).
\end{itemize} 
\end{definition}
For any fixed $T>0$, $p\geq 1$, $q\geq 1$, $s\geq 0$ and for any fixed (vector-valued) function $f: [0,T]\times \R^2\mapsto \R^m$, $m\geq 1$, we denote
\begin{align}\label{norm}
\|f\|_{L^p_T X_x}:=\big\|  \|f(t)\|_{X_x(\R^2; \R^m)}\big\|_{L^p_t([0,T])}
\hbox{ with }X= H^s\hbox{ or }L^q.
\end{align}
The functional space $L^p([0,T];H^s(\R^2; \R^m))$ consists of all functions $f:[0,\infty)\times \R^2\to\R^m$ satisfying $\|f\|_{L^p_TH^s_x}<\infty$.
We have the following existence, uniqueness   as well as   global-in-time regularity results for the solutions of the Cauchy problem for the Boussinesq equations \eqref{BS}-\eqref{BS-IC}-\eqref{ab} on the whole two dimensional space $\R^2$. 
\begin{theorem}[Existence, uniqueness \& Global-in-time regularity]\label{thm:well-posedness}
For any initial data $\theta_0\in L^2(\R^2)$ and $u_0\in (L^2(\R^2))^2$, there exists a  global-in-time weak solution 
$$(\theta,u)\in C([0,\infty); (L^2(\R^2))^3)\cap L^2_\loc([0,\infty); (H^1(\R^2))^3)$$ 
of the initial value problem  \eqref{BS}-\eqref{BS-IC}-\eqref{ab}.

If $\theta_0\in H^1(\R^2)$, $u_0\in (L^2(\R^2))^2$ and the functions $a,b\in C^2_b(\R)$ have finite first and second derivatives, then the weak solution  is indeed unique, and satisfies
\begin{align*}
  &\theta\in C([0,\infty);H^1(\R^2))\cap L^2_\loc([0,\infty);H^2(\R^2)), 
\end{align*}
as well as  the following energy estimates for any given $T>0$, 
\begin{align}
&\|u\|^2_{L^\infty_TL^2_x}+\|\nabla u\|^2_{L^2_T L^2_x}\le C\big(T\|\theta_0\|^2_{L^2}+\|u_0\|^2_{L^2}\big),\label{energy:u:1}\\
&\|\theta\|_{L^\infty_TH^1_x}^2+
   \|(\d_t \theta, \nabla^2\theta)\|_{L^2_TL^2_x}^2
  \notag \\
&   \le C\|\theta_0\|_{H^1}^2(1+\|\nabla \theta_0\|_{L^2}^2)\exp\big(C(T^2\|\theta_0\|^4_{L^2}+\|u_0\|^4_{L^2})\big),\label{energy:theta:1}
\end{align}
where $C$ is a positive constant depending only on $\|a\|_\Lip, \kappa_\ast, \kappa^\ast, \mu_\ast$.


Furthermore, the general $H^s$-regularities can   be propagated globally in time in the following sense:  For any initial data  (see the grey  closed unbounded quadrangle in Figure \ref{fig:s} for the admissible regularity exponent range)
\begin{equation}\label{s}\begin{split}
&(\theta_0, u_0)\in H^{s_\theta}(\R^2)\times (H^{s_u}(\R^2))^2
\\
&\hbox{ with }(s_\theta, s_u)\in \{(s_\theta, s_u)\subset [1,\infty)\times [0,\infty)\,|\, s_u-1\leq s_\theta\leq s_u+2\}
\end{split}\end{equation} 
and the functions $a\in C^2_b\cap C^{[s_\theta]+1}$, $b\in C^2_b\cap C^{[s_u]+1}$, the unique solution $(\theta, u)$ stays in 
\begin{equation}\label{regularity}
 C([0,\infty);H^{s_\theta}(\R^2)\times(H^{s_u}(\R^2))^2)\cap L^2_\loc([0,\infty);H^{s_\theta+1}(\R^2)\times (H^{s_u+1}(\R^2))^2).
\end{equation}
\end{theorem} 

 Theorem \ref{thm:well-posedness} will be proved in Section \ref{sec:proof}.
The proof of the existence of weak solutions is rather standard, and we are going to sketch the proof in Subsection \ref{subs:existence} for the reason of completeness, as we did not find the proof in the literature.
As mentioned before, some well-posedness results have already been established for smooth data in the bounded domain case (see \eqref{LB} above in e.g. \cite{DiazGaliano, Duvaut:Lions, Goncharova, Lorca:Boldrini}) or in  smoother functional frameworks in the whole space case (see \eqref{WZ} above in e.g. \cite{Wang:Zhang}).  
We are going to focus on the proofs of the uniqueness result and the global-in-time regularity in the low regularity regimes in Subsection \ref{subs:uniqueness} and Subsection \ref{sec:fractional} respectively, where different regularity exponents for different unknowns are permitted.
The commutator estimates as well as the composition estimates in Lemma \ref{Liao} will play an important role, and the a priori estimates for a general linear parabolic equation in Lemma \ref{lem:parabolic} will be of independent interest.

 We conclude this introduction part with  several remarks on the results in Theorem \ref{thm:well-posedness}.
\begin{remark}[Optimality of the regularity exponent ranges in Theorem \ref{thm:well-posedness}]\label{rem:optimality}
We are going to follow the standard procedure to show the existence of weak solutions for $L^2$-initial data  by use of the a priori energy (in)equalities \eqref{EE:theta} and \eqref{EE:u} (see Subsection \ref{subs:existence}).

If we take the difference between two different weak solutions $(\theta_1, u_1)$ and $(\theta_2, u_2)$, the difference of the nonlinear viscosity term $\div(\mu Su)$ in the $u$-equation will become
$$
\div((\mu_1-\mu_2) Su_1)+\div(\mu_2 S(u_1-u_2)),
$$
which stays in $L^2_\loc([0,\infty); (H^{-1}(\R^2))^2)$ provided with
\begin{align*}
&u_1\hbox{ and }(u_1-u_2)\in L^2_\loc([0,\infty); (H^{1}(\R^2))^2), 
\\
&\hbox{and }\mu_1-\mu_2\in L^\infty_\loc\bigl([0,\infty); H^1(\R^2)\bigr)\subset L^\infty_\loc\bigl ([0,\infty); (L^1(\R^2))'=\hbox{BMO}(\R^2)\bigr).
\end{align*}  
Therefore in order to ensure the $L^2_x$-Estimate for the velocity difference $(u_1-u_2)$, we require the $H^1_x$-Estimate for the  temperature difference in $(\mu_1-\mu_2)$.
And hence the initial condition $\theta_0\in H^1_x$ (i.e. $s_\theta\geq 1$ above) is required for the proof of the uniqueness result   (see Subsection \ref{subs:uniqueness}). 

Under the lower-regularity assumption $\theta_0\in H^s_x$ with $0<s<1$,  the coefficients $\kappa, \mu$ are not expected to be continuous uniformly in time, and hence no uniqueness or  $H^s$-regularity results for $\theta$ or $H^{s_1}$, $s_1>0$-regularity results for $u$  are expected. 
{}
Nevertheless with constant diffusion coefficients (e.g. $\kappa=\mu=1$), the uniqueness result for the weak solutions holds true by virtue of the $L^2_x$-energy (in)equalities (similar as the classical global-in-time well-posedness result for the classical two dimensional incompressible Navier-Stokes equations).
{}
Furthermore, if $\kappa=1$ is a positive constant, then the $H^s_x$, $s\in (0,1)$-Estimate for $\theta$ holds true, provided with $u\in L^4_\loc(L^4_x(\R^2))^2$ (or with $u_0\in (L^2(\R^2))^2$), simply by an interpolation argument between \eqref{EE:theta} and \eqref{energy:theta:1}.
{}
Similarly if $\mu=1$ is a positive constant, then  the $H^s_x$, $s>0$-Estimate for $u$ holds true, provided with $\theta\in L^2_\loc(H^{s-1}_x(\R^2))$.
{}
Thus with constant diffusion coefficients (e.g. $\kappa=\mu=1$), the Sobolev regularities 
$$
(\theta_0, u_0)\in (H^s(\R^2))\times(L^2(\R^2))^2\hbox{ or }(L^2(\R^2))\times  (H^s(\R^2))^2,\quad 0<s\leq 1
$$
can be propagated globally in time, and the admissible regularity exponent set \eqref{s} extends itself indeed to the   closed set consisting of all non-negative admissible regularity exponents:
\begin{equation*}\begin{split}
 (s_\theta, s_u)\in \{(s_\theta, s_u)\subset [0,\infty)\times [0,\infty)\,|\, s_u-1\leq s_\theta\leq s_u+2\}.
\end{split}\end{equation*} 

In order to propagate the $H^{s_\theta}$, $s_\theta\geq 2$-regularity of $\theta$, we require the  transport term $u\cdot\nabla\theta$ in the $\theta$-equation to be at least in $L^2_\loc([0,\infty); H^{s_\theta-1}_x)$, which requires  $u\in L^2_\loc([0,\infty); H^{s_\theta-1}_x)$ and hence the initial assumption $u_0\in H^{s_u}$ with the restriction  $s_u\geq s_\theta-2$ (as there is a gain of regularity of oder $1$ when taking $L^2$-norm in the time variable in general).
{}
Similarly, in order the propagate the $H^{s_u}$, $s_u\geq 2$-regularity of $u$, we require the viscosity term $\div(\mu Su)$ in the $u$-equation to be at least in  $L^2_\loc([0,\infty); H^{s_u-1}_x)$, which requires $\mu Su\in L^2_\loc([0,\infty); H^{s_u}_x)$
and hence   the initial assumption  $\theta_0\in H^{s_\theta}$ with the restriction $s_\theta\geq s_u-1$.
\end{remark}

\begin{remark}[Precise $H^s_x$-Estimates in the high regularity regime]\label{rem:Hs}
The global-in-time regularity in the high regularity regime \eqref{s}-\eqref{regularity} follows immediately from  the following  \emph{borderline} a priori estimates:
\begin{itemize}
\item
If  $\theta_0\in H^{s}(\R^2)$, $u_0\in (L^2(\R^2))^2$ with $s\in (1,2]$ and the function $a\in C^2_b(\R)$, then 
for $s\in (1,2)$ it holds
\begin{equation}\label{theta:Hsnu}
\begin{aligned}
   & \|\theta\|^2_{ L^\infty_T H^{s}_x }+ \|\nabla\theta\|^2_{ L^2_T H^{s}_x }\le
    C(\kappa_\ast) \|\theta_0\|_{H^s_x}^2\times
   \\
   &\times \exp\Bigl(C(\kappa_\ast, s,   \|a\|_{C^2}, \|\theta\|_{L^\infty_TH^1_x})\bigl(\|\nabla u\|_{L^2_T L^2_x}^2+ \|\nabla\theta\|_{L^2_T H^1_x}^2\bigr)\Bigr),
\end{aligned}   
\end{equation}
 and for $s=2$ it holds
\begin{equation}\label{theta:Hs2}
\begin{aligned}
&\|\theta\|_{L^\infty_T H^2_x}^2+\|\nabla\theta\|_{L^2_T H^2_x}^2
\le C(\kappa_\ast,\|a\|_{C^2}, \kappa^\ast) \|\theta_0\|^2_{H^2}(1+\|\nabla\theta_0\|^2_{L^2})^2
\\
&\qquad  \times
 \exp\Bigl(C(\kappa_\ast,\|a\|_{\Lip})(\|\nabla u\|_{L^2_T L^2_x}^2+\|u\|_{L^4_T L^4_x}^4+\|\nabla\theta \|_{L^4_T L^4_x}^4)\Bigr) .
\end{aligned}   
\end{equation}

\item
If $\theta_0\in H^1(\R^2)$, $u_0\in (H^{s}(\R^2))^2$ with $s\in (0,2]$ and the function $b\in C^2_b(\R)$,  
then  for $s\in (0,2)$ it holds
  \begin{equation}\label{u:Hsnu}
\begin{aligned} 
&\|u\|_{L^\infty_T H^s_x}^2+\|\nabla u\|_{L^2_T H^{s}_x}^2
\leq  C(\mu_\ast) ( \|u_0\|_{H^s_x}^2+T\|\theta_0\|_{  L^2_x}^2+\|\theta\|_{L^2_T H^{s-1}_x}^2)
\\
&\times \exp\Bigl(C(\mu_\ast,s, \|b\|_{C^2}, \|\theta\|_{L^\infty_T H^1_x})\bigl(\|\nabla u\|_{L^2_T L^2_x}^2+ \|\nabla\theta\|_{L^2_T H^1_x}^2\bigr)\Bigr),
\end{aligned}   
\end{equation}  
  and for $s=2$ it holds
\begin{equation}\label{u:Hs2}\begin{split}
&\| u\|_{L^\infty_T H^2_x}^2+\| \nabla  u\|_{L^2_T H^2_x}^2 
 \leq  (\| u\|_{L^\infty_T H^1_x}^2+\| \nabla  u\|_{L^2_T H^1_x}^2)
 +C(\mu_\ast, \|b\|_{C^2})\times
 \\
 &
(1+ \|\Delta u_0\|_{L^2_x}^2) (1+\| u\|_{L^\infty_T H^1_x}^2+\| \nabla  u\|_{L^2_T H^1_x}^2)^2
(1+\| \theta\|_{L^\infty_T H^1_x}^2+\| \nabla  \theta\|_{L^2_T H^1_x}^2)  
 \\
&\times  \exp\Bigl(C(\mu_\ast, \|b\|_{C^2})(\|(u, \nabla \theta)\|_{L^4_T L^4_x}^4+\|\nabla^2\theta \|_{L^2_T L^2_x}^2)\Bigr).
\end{split}\end{equation}


\item If  $\theta_0\in H^{s}(\R^2)$, $u_0\in (H^{s-2}(\R^2))^2$ with $s>2$ and the function $a\in  C^{[s]+1}$,   then 
for $s\in (2,3)$ it holds
\begin{equation}\label{theta:Hs23}\begin{split}
&\|\theta\|_{L^\infty_T H^s_x}^2+\|\nabla\theta\|_{L^2_T H^{s}_x}^2
\leq  C(\kappa_\ast )\|\theta_0\|_{H^s_x}^2\times 
\\
&\quad
\times \exp\Bigl(C(\kappa_\ast,s, a, \|\theta\|_{L^\infty_T L^\infty_x})(\| u\|_{L^2_T H^{s-1}_x}^2+ \|\nabla\theta\|_{L^2_T L^\infty_x}^2)\Bigr),
\end{split}\end{equation}
  and  for $s\geq 3$ it holds
\begin{equation}\label{theta:Hs3}\begin{split}
&\|\theta\|_{L^\infty_T H^s_x}^2+\|\nabla\theta\|_{L^2_T H^{s}_x}^2
 \leq  
C(\kappa_\ast,s)  (\|\theta_0\|_{H^s_x}^2+\|\nabla\theta\|_{L^\infty_T L^\infty_x}^2\|\nabla u\|_{L^2_T H^{s-2}_x}^2)
\\
&\qquad  
\times\exp(C(\kappa_\ast,s, a, \|\theta\|_{L^\infty_T L^\infty_x})(\|\nabla u\|_{L^2_T L^\infty_x}^2+ \|\nabla\theta\|_{L^2_T L^\infty_x}^2)). 
\end{split}\end{equation} 
 
\item If  $\theta_0\in H^{s-1}(\R^2)$, $u_0\in (H^{s}(\R^2))^2$ with  $s>2$ and the function  $b\in C^{[s]+1}$, then 
for $s\in (2,3)$ it holds  
\begin{equation}\label{u:Hs23}\begin{split}
&\|u\|_{L^\infty_T H^s_x}^2+\|\nabla u\|_{L^2_T H^{s}_x}^2
\leq  C(\mu_\ast)  ( \|u_0\|_{H^s_x}^2+T\|\theta\|_{L^\infty_T H^{s-1}_x}^2)\\
& \times\exp(C(\mu_\ast, s, b, \|\theta\|_{L^\infty_TL^\infty_x})(\|\nabla u\|_{L^2_T H^1_x}^2+ \|\nabla\theta\|_{L^2_T H^{s-1}_x}^2)),
\end{split}\end{equation}
and for $s\geq 3$ it holds
\begin{equation}\label{u:Hs3}\begin{split}
&\|u\|_{L^\infty_T H^s_x}^2+\|\nabla u\|_{L^2_T H^{s}_x}^2
\\
&\leq  
C(\mu_\ast ) (\|u_0\|_{H^s_x}^2+T\|\theta\|_{L^\infty_T H^{s-1}_x}^2
 +\|\nabla u\|_{L^\infty_T L^\infty_x}^2\|\nabla \theta\|_{L^2_T H^{s-1}_x}^2) 
\times \\
& \times\exp( C(\mu_\ast, s, b, \|\theta\|_{L^\infty_TL^\infty_x})(\|\nabla u\|_{L^2_T L^\infty_x}^2+ \|\nabla\theta\|_{L^2_T L^\infty_x}^2)). 
\end{split}\end{equation}  
\end{itemize}
\end{remark} 


\begin{remark}[$L^2$-in time Estimate VS $L^1$-in time Estimate]\label{rem:L2}
   Instead of the classical $L^\infty_t H^s_x\cap L^1_t H^{s+2}_x$-type estimate in the literature, we derive $L^\infty_t H^s_x\cap L^2_t H^{s+1}_x$-type estimate here, since  e.g. only the $L^2_t \dot H^1_x$-a priori estimate for the velocity vector field is available from  the energy estimates (roughly speaking, the $L^2_t$-in time norm asks less spacial regularity on the coefficients).
    See Lemma \ref{lem:parabolic} below for the a priori $H^s_x$, $s\in (0,2)$-estimates for a general linear parabolic equation with divergence-free $L^2_t H^1_x$-velocity vector field, which is of independent interest.  
   
It is in general not true that $\theta\in L^1_t H^{s+2}_x$ (or $u\in L^1_t H^{s+2}_x$) in the low regularity regime, although it holds straightforward in the high regularity regime.
\end{remark} 

\begin{remark}[Remarks on the smoothness assumptions on the functions $a,b$]
It is common to assume smooth heat conductivity law and viscosity law  \cite[I]{Perez} in modelling fluids.

The Lipschitz continuity assumption $a,b\in \Lip$ is enough for the existence result as well as for the $H^1\times L^2$-Estimates \eqref{energy:u:1}-\eqref{energy:theta:1} in Theorem \ref{thm:well-posedness}.

As for the uniqueness result, due to the following $\dot H^1_x$-Estimate for the difference of the diffusion coefficinets
\begin{align*}
\|a(\theta_1)-a(\theta_2)\|_{\dot H^1_x}
&\leq \|(a'(\theta_1)-a'(\theta_2))\nabla\theta_1\|_{L^2_x}+\|a(\theta_2)\nabla(\theta_1-\theta_2)\|_{L^2_x}
\\
&\leq 
 \|a'\|_{\Lip} \|\nabla\theta_1\|_{L^4_x} \|\theta_1-\theta_2\|_{L^4_x}+\|a\|_{L^\infty}  \| \nabla(\theta_1-\theta_2)\|_{L^2_x},
\end{align*}
the Lipschitz continuity assumptions $a',b'\in \Lip$ are required.

The  dependance on the function $a$ of the constants $C$ in  \eqref{theta:Hs23}-\eqref{theta:Hs3} reads precisely as (similarly for the constants in \eqref{u:Hs23}-\eqref{u:Hs3})
$$
\sup_{k=0,\cdots, [s]+1}\,\sup_{|y|\leq c\|\theta\|_{L^\infty_T L^\infty_x}} \bigl| \frac{d}{dy^k} a(y)\bigr|.
$$

For the integer regularity exponents, we can simply derive the energy estimates by integration by parts (instead of the application of the commutator estimates or the composition estimates in Lemma \ref{Liao} below), such that the requirement for $a\in C^{[s_\theta]+1}$ and $b\in C^{[s_u]+1}$ can be relaxed, see e.g. \eqref{theta:Hs2}, \eqref{u:Hs2}. 
\end{remark}

\section{Proofs}\label{sec:proof}
Recall the Cauchy problem for the two dimensional Boussinesq equations \eqref{BS}-\eqref{ab}
\begin{equation}\label{Boussinesq}
\left\{
\begin{aligned}
&\d_t\theta+ u\cdot \nabla \theta-\di(\kappa \nabla \theta)=0,\\
&\d_t u+u\cdot \nabla u-\di(\mu  Su)+\nabla\Pi= \theta \vec{e_2},\\
&\di u=0,\\
&(\theta,u)\mid_{t=0}=(\theta_0,u_0) ,
\end{aligned}
\right.
\end{equation}
where $\kappa=a(\theta)\in C^1_b (\R; [\kappa_\ast, \kappa^\ast])$, $\mu=b(\theta)\in C^1_b(\R; [\mu_\ast, \mu^\ast])$ with   $\kappa_\ast, \kappa^\ast, \mu_\ast, \mu^\ast$ four positive constants.

We are going to show the existence, uniqueness as well as the global-in-time regularity results in Theorem \ref{thm:well-posedness} in Subsection \ref{subs:existence}, Subsection \ref{subs:uniqueness} and Subsection \ref{sec:fractional} respectively.

Recall the definition of the  $\|\cdot\|_{L^q_T X_x}$-norm in \eqref{norm}.
The Gagliardo-Nirenberg's inequality
\begin{equation}\label{GN}
\lVert f\rVert_{L^4_TL^4_x(\R^2)}\leq C\lVert f\rVert_{L^2_TL^2_x(\R^2)}^{\frac12}\lVert \nabla f\rVert_{L^2_TL^2_x(\R^2)}^{\frac12}
\end{equation}
as well as the equivalence relations between   the norms   
\begin{equation}\label{equiv:norm}\begin{split}
&\|Su\|_{L^2_x(\R^2)}^2=2\|\nabla u\|_{L^2_x(\R^2)}^2\hbox{ if }\div u=0,
\\
&\lVert \Delta\eta\rVert_{L^2_x(\R^2)}\sim \lVert \nabla^2\eta\rVert_{L^2_x(\R^2)}
\end{split}\end{equation}
will be used freely in the proof.

\subsection{Existence of weak solutions if $(\theta_0,u_0) \in   (L^2(\R^2))^3$}\label{subs:existence} 
We will follow the standard procedure to show the existence   of the weak solutions under the initial condition
$$(\theta_0,u_0) \in L^2(\R^2)\times (L^2(\R^2))^2,$$
 namely
\begin{enumerate}[Step 1]
\item We construct a sequence of approximate solutions, which satisfy the energy estimate uniformly. 
\item We show the convergence of this approximate solution sequence to a weak solution and study the property of the weak solution.
\end{enumerate}
We are going to sketch the proof and  pay attention to the low-regulartiy assumptions.

\begin{paragraph}{Step 1: Construction of approximate solutions with uniform bounds}
We use the Friedrich's method to construct a sequence of approximate solutions. We consider the following system of $(\theta_n,u_n)$
\begin{equation}\label{app2}
\left\{
\begin{aligned}
&\d_t\theta_n+ P_n( u_n
\cdot \nabla \theta
_n)-P_n\di(\kappa_n \nabla \theta_n)=0,\\
&\d_t u_n+P_n \PP ( u_n\cdot \nabla u_n)-P_n \PP \di(\mu_n Su_n)= \PP( \theta_n \vec{e_2}),\\
&u_n(0,x)=P_nu_0(x),\quad \theta_n(0,x)=P_n\theta_0(x),
\end{aligned}
\right.
\end{equation}
where $\kappa_n=a(\theta_n)$ and $\mu_n=b(\theta_n)$. The operator $P_n$, $n\in\N$, is the low-frequency cut-off operator which is defined as follows
$$P_nf(x)=\mathcal F^{-1}(\mathbb{1}_{B_n}(\xi)\mathcal{F} f(\xi))(x),$$
where $B_n\subset\R^2$ is the disk with center at $0$ and radius $n$, and $\mathcal F, \mathcal F^{-1}$ are the standard Fourier and inverse Fourier transformation.
The operator $\PP$ in \eqref{app2} denotes the Leray-Helmholtz projector
on the $\R^2$, which  decomposes the tempered distributions $v\in \cS'(\R^2; \R^2)$ into div-free and curl-free parts as follows
\begin{equation}\label{LerayProjector}
v=\nabla^\perp V_1+\nabla V_2,
\end{equation}
where 
$$\nabla^\perp V_1=-\nabla^\perp(-\Delta)^{-1}\nabla^\perp\cdot v=:\PP v,
\quad \nabla V_2=-\nabla(-\Delta)^{-1}\nabla\cdot v= (1-\PP)v$$
 with $\nabla^\perp=(\d_{x_2}, \, -\d_{x_1})^T$. 
Notice that   $\PP$ maps $L^p(\R^2; \R^2)$ into itself  for any $p\in (1,\infty)$ and it is commutative with the projection operator $P_n$.

 We define the Banach spaces $L^2_n$ and $L^{2,\sigma}_n$ as following
\begin{align*}
L^2_n(\R^2)&=\{f\in L^2(\R^2)\mid f=P_nf\},\\
L^{2,\sigma}_n(\R^2)&=\{f\in (L^2_n(\R^2))^2\mid \div_x(f)=0\}.
\end{align*} 
The system \eqref{app2} turns out to be an ordinary differential equation system in $L^2_n(\R^2)\times (L^{2,\sigma}_n(\R^2))^2 $. Indeed, the following estimates hold  \begin{align*}
&   \|P_n( u_n \cdot \nabla \theta
_n)-P_n\di(\kappa_n \nabla \theta_n)\|_{L^2_x} \\
&\qquad \le C(\|a\|_{\Lip}) n^4 (\|u_n\|_{L^2_x}+\|\theta_n\|_{L^2_x}+\kappa^\ast)\|\theta_n\|_{L^2_x},\\
& \|P_n \PP ( u_n\cdot \nabla u_n)-P_n \PP \di(\mu_n Su_n)\|_{L^2_x} \\
&\qquad
 \le C(\|b\|_{\Lip}) n^4 (\|u_n\|_{L^2_x}+\|\theta_n\|_{L^2_x}+\mu^\ast)\|u_n\|_{L^2_x}.
\end{align*}
Hence, for any $n\in \N$, there exists $T_n>0$ such that the system \eqref{app2} has a solution $(\theta_n,u_n)\in C^\infty([0,T_n];L_n^2(\R^2))\times(C^\infty([0,T_n];L_n^{2,\sigma}(\R^2)))^2$.

We take the $L^2$ inner product of the equation $\eqref{app2}_1$ and $\theta_n$ to derive
$$\frac12\frac{d}{dt}\int_{\R^2} \theta_n^2+\int_{\R^2}\kappa_n |\nabla \theta_n|^2=0.$$
Then the following uniform estimate for $(\theta_n)$ holds
\begin{equation}\label{uniform:theta}
\frac12\|\theta_n\|_{L^\infty_T L^2_x}^2+\kappa_\ast \|\nabla \theta_n\|_{L^2_TL^2_x}^2\,dt
\leq \frac12 \|P_n\theta_0\|_{L^2_x}^2\le \frac12 \|\theta_0\|_{L^2_x}^2,\quad \forall T>0.
\end{equation}
Similarly we take the $L^2$ inner product of the equation $\eqref{app2}_2$ and $u_n$ to derive
\begin{align*}
\frac12  \frac{d}{dt}\|u_n\|_{L^2_x}^2+\frac12\|  \mu_n S u_n\|^2_{L^2_x} &\le  \|\theta_n\|_{L^2_x}\|u_n\|_{L^2_x}\le \frac12(T\| \theta_n\|_{L^2_x}^2+\frac1T\| u_n\|_{L^2_x}^{2}), 
\end{align*}
for all positive times $T>0$,
and thus by Gronwall's inequality we arrive at the following uniform estimate for $(u_n)$ (noticing $\|Su_n\|_{L^2_x}^2=2\|\nabla u_n\|_{L^2_x}^2$)
\begin{equation}\label{uniform:u}
\frac12 \|u_n\|_{L^\infty_TL^2_x}^2+\mu_\ast \|\nabla u_n\|_{L^2_T L^2_x}^2
\le  \frac e2  (   T\|\theta_0\|_{L^2_x}^2+\|u_0\|_{L^2_x}^2),\quad \forall T>0.   
\end{equation}
Thus the approximate solutions $(\theta_n, u_n)$ exist for all positive times.
\end{paragraph}

\begin{paragraph}{Step 2: Passing to the limit}
By the above uniform bounds \eqref{uniform:theta}-\eqref{uniform:u} there exists a subsequence, still denote by $(\theta_{n},u_{n})$,
converging weakly to a limit $(\theta, u)\in L^\infty_\loc([0,\infty); (L^2_x)^3)\cap L^2_\loc([0,\infty); ( H^1_x)^3)$:
\begin{align*}
    &\theta_n\overset{*}{\rightharpoonup} \theta \quad \text{in }  L^\infty_\loc([0,\infty); L^2(\R^2)),\quad \nabla \theta_n\rightharpoonup \nabla\theta \quad\text{in } L^2_\loc([0,\infty); L^2(\R^2)),\\
     &u_n\overset{*}{\rightharpoonup} u \quad \text{in } L^\infty_\loc([0,\infty); (L^2(\R^2))^2),
      \quad \nabla u_n\rightharpoonup \nabla u \quad\text{in } L^2_\loc([0,\infty);  (L^2(\R^2))^2).
\end{align*}
Since by the Gagliardo-Nirenberg's inequality $(\theta_n, u_n)$ is a bounded sequence in $L^4_TL^4_x$ for any $T>0$, the sequence of the time derivatives   $(\d_t\theta_n, \d_t u_n)$ is   bounded in $L^2_T(H^{-1}_x)$ (by use of the equations in \eqref{app2}), and hence $\{(\theta_n, u_n)\}$ is relatively compact in $L^p_T L^2_x(B_R)$ for any fixed disk $B_R\subset\R^2$ for all $p\in [1,\infty)$, which implies the pointwise convergence (up to a subsequence) 
$$\theta_n\rightarrow \theta,\quad u_n\rightarrow u \hbox{ for almost every }t\in \R^+,\quad x\in \R^2,$$
as well as the convergence of the nonlinear terms (noticing e.g. $u_n\varphi\rightarrow u\varphi$ in $L^4_TL^4_x$ for fixed $\varphi\in C_c^\infty((0,T)\times\R^2)$)
\begin{align*}
&u_n\theta_n\rightarrow u\theta,\quad u_n\otimes u_n\rightarrow u\otimes u\hbox{ in }\cD'((0,T)\times\R^2)
\hbox{ and hence weakly in }L^2_TL^2_x.
\end{align*} 
Consequently, $\kappa_n=a(\theta_n)\rightarrow \kappa=a(\theta)$ and $\mu_n=b(\theta_n)\rightarrow \mu=b(\theta)$ almost everywhere and 
\begin{align*} 
&\kappa_n\nabla\theta_n\rightharpoonup \kappa\nabla \theta,\quad \mu_n S u_n\rightharpoonup \mu S u
\hbox{ in }L^2_TL^2_x.
\end{align*}

Thus the equation (noticing $P_n\rightarrow \hbox{Id}$ as an operator from $H^s(\R^2)$  to itself)
$$\d_t\theta+\div(u\theta)-\div(\kappa\nabla\theta)=0\hbox{ holds in  }L^2_\loc((0,\infty); H^{-1}_x(\R^2)),$$
and we can test it by $\theta\in L^2_\loc([0,\infty); H^1_x)$ to arrive at the energy equality \eqref{EE:theta} for $\theta$, such that $\theta|_{t=0}=\theta_0$ and $\theta\in C([0,\infty); L^2_x)$ hold true.

Similarly, the equation
   \begin{align}\label{eq:u}
   \d_t u+\PP\div(u\otimes u-\mu Su)=\PP(\theta \vec{e}_2) \hbox{ holds in }L^2_\loc((0,\infty); (H^{-1}_x(\R^2))^2),
   \end{align}
   and we can test it by the divergence-free velocity field $u\in L^2_\loc((0,\infty); (H^{1}_x(\R^2))^2)$ to arrive at the energy equality \eqref{EE:u}, which implies $u\in C([0,\infty); (L^2_x(\R^2))^2)$ and $u|_{t=0}=u_0$.
   We take the solution $\Pi$ of the Poisson equation
   \begin{equation}\label{eq:Pi}
   \Delta\Pi=\div(1-\PP)(\theta\vec{e}_2-\div(u\otimes u-\mu Su))
   \end{equation} 
   under the renormalisation condition $\int_{B_1}\Pi \dx=0$, such that
   $$
   \nabla\Pi= (1-\PP)(\theta\vec{e}_2-\div(u\otimes u-\mu Su))
   \in L^2_\loc((0,\infty); (H^{-1}_x(\R^2))^2),
   $$
   and the equation \eqref{BS:u} holds in $L^2_\loc((0,\infty); (H^{-1}_x(\R^2))^2)$.
   \end{paragraph}

\subsection{Energy estimates   \& Uniqueness of the weak solutions if $(\theta_0,u_0) \in H^1(\R^2)\times (L^2(\R^2))^2$}\label{subs:uniqueness} 
We first introduce a scalar function $\eta$, which  is  given in terms of the temperature function as follows (recalling  $\kappa=a(\theta)\in C^1_b(\R; [\kappa_\ast, \kappa^\ast])$)
\begin{equation}\label{eta}
\eta=A(\theta),\hbox{ with } A(z):=\int_0^z a(\alpha)\,d\alpha \hbox{ the primitive function of }a.
\end{equation}
As $A'(\theta)=a(\theta)\geq\kappa_\ast>0$, the function $A$ is invertible and we can write 
\begin{equation}\label{eta,theta}
\theta=A^{-1}(\eta),
\end{equation}
where $(A^{-1})'(\eta)=\frac{1}{a(A^{-1}(\eta))}\leq \frac1{\kappa^\ast}$.
We have  the following equivalence relations \footnote{We can easily compute  
\begin{align*}
&\nabla\eta=a(\theta)\nabla\theta,\quad  \nabla\theta=\frac{1}{a(A^{-1}(\eta))}\nabla\eta,\\
&\nabla^2\eta=a'(\theta)\nabla\theta\otimes\nabla\theta+a(\theta)\nabla^2\theta,
\quad
\nabla^2\theta=-\frac{a'(A^{-1}(\eta))}{a^3(A^{-1}(\eta))}\nabla\eta\otimes\nabla\eta+\frac{1}{a(A^{-1}(\eta))}\nabla^2\eta.
\end{align*} }
\begin{equation}\label{eq:eta,theta}\begin{split}
&\kappa_\ast\|\theta\|_{L^2_x}\leq \|\eta\|_{L^2_x}\leq\kappa^\ast\|\theta\|_{L^2_x},\\
&\kappa_\ast\|\nabla\theta\|_{L^2_x}\leq \|\nabla \eta\|_{L^2_x}=\|a(\theta)\nabla \theta\|_{L^2_x}\leq \kappa^\ast\|\nabla\theta\|_{L^2_x},
\\
& \kappa_\ast\|\d_t\theta\|_{L^2_x}\leq \|\d_t \eta\|_{L^2_x}=\|a(\theta)\d_t \theta\|_{L^2_x}\leq \kappa^\ast\|\d_t\theta\|_{L^2_x},
\\
&\|\nabla^2\eta\|_{L^2_x}\leq  \|a\|_{\Lip}\|\nabla\theta\|_{L^4_x}^2+\kappa^\ast\|\nabla^2 \theta\|_{L^2_x}
\leq  (C\|a\|_{\Lip}\|\nabla\theta\|_{L^2_x}+\kappa^\ast)\|\nabla^2 \theta\|_{L^2_x},
\\
& \|\nabla^2 \theta\|_{L^2_x}\leq \frac{\|a\|_{\Lip}}{\kappa_\ast^3}\|\nabla\eta\|_{L^4_x}^2+\frac{1}{\kappa_\ast}\|\nabla^2 \eta\|_{L^2_x}\leq (C\frac{\|a\|_{\Lip}}{\kappa_\ast^3}\|\nabla\eta\|_{L^2_x}+\frac{1}{\kappa_\ast})\|\nabla^2 \eta\|_{L^2_x}.
\end{split}\end{equation}   
That is,
\begin{equation}\label{equiv:eta,theta}
\theta(t,\cdot)\in H^k_x(\R^2)\Leftrightarrow \eta(t,\cdot)\in H^k_x(\R^2),\quad k=0,1,2.
\end{equation}

Let $(\theta, u)\in C([0,\infty); (L^2(\R^2))^3)\cap L^2_\loc([0,\infty); (H^1(\R^2))^3)$ be a weak solution of the Cauchy problem \eqref{Boussinesq} in the sense of Definition \ref{def:weak} with
\begin{equation}\label{theta:eq}
\d_t \theta+u\cdot\nabla\theta-\div(\kappa\nabla\theta)=0
\hbox{ holding in } L^2_\loc([0,\infty); H^{-1}_x(\R^2)).
\end{equation}
Since   $Y:=L^\infty_{t,x}([0,\infty)\times\R^2)\cap  L^2_\loc([0,\infty); H^1_x(\R^2))$ is an algebra (in the sense that the product of any two elements in $Y$ still belongs to $Y$),    we can multiply the above $\theta$-equation by   $\kappa=a(\theta)$ (with  $ a(\theta)-a(0)\in Y$), to  arrive   at the parabolic equation for $\eta=A(\theta)\in C([0,\infty); L^2(\R^2))\cap L^2_\loc([0,\infty); H^1(\R^2))$:
\begin{equation}\label{eta:eq}
 \d_t \eta+u\cdot\nabla \eta -\kappa\Delta \eta=0\hbox{ holding in the dual space }Y'.
 \end{equation}

 We are going to derive the   $H^1$-Estimate for $\eta$ (and hence for $\theta$) as well as the  $L^2$-Estimate for $u$ first.
Then we will show the uniqueness result of the weak solutions by considering the difference of two possible weak solutions.
 The procedure is standard (see  e.g. Section 2 \cite{Liao}) 
  and we are going to sketch the proof.

\begin{paragraph}{$H^1\times L^2$-Estimate for $(\theta, u)$} 
By virtue of the energy equalities \eqref{EE:theta} and \eqref{EE:u} and the derivation of the uniform estimates \eqref{uniform:theta} and \eqref{uniform:u}, we have the $L^2$-Estimate
\begin{equation}\label{theta:LinftyL2}
\|\theta\|^2_{L^\infty_TL^2_x}+\|\nabla \theta\|^2_{L^2_T L^2_x}\le C(\kappa_\ast) \|\theta_0\|^2_{L^2},
\end{equation}
and the $L^2$-Estimate
 \eqref{energy:u:1} for $u$.
By Gagliardo-Nirenberg's inequality \eqref{GN} it holds 
\begin{equation}\label{u:L4}
\|u\|_{L^4_TL^4_x} \le C(\mu_\ast)\big(\sqrt T\|\theta_0\|_{L^2}+\|u_0\|_{L^2}\big).
\end{equation}

We assume a priori that the function $\eta$ is smooth and decay sufficiently fast at infinity.  
We test the $\eta$-equation \eqref{eta:eq} by $\Delta\eta$ to derive by integration by parts that
\begin{align*}
\frac12\frac{d}{dt}\int_{\R^2} |\nabla\eta|^2\dx+\int_{\R^2}\kappa|\Delta\eta|^2\dx
&=\int_{\R^2}u\cdot\nabla\eta\Delta\eta\dx
\leq \|u\|_{L^4_x}\|\nabla\eta\|_{L^4_x}\|\Delta\eta\|_{L^2_x}.
\end{align*}
By Gagliardo-Nirenberg's inequality \eqref{GN}, the equivalence $\lVert \Delta\eta\rVert_{L^2_x(\R^2)}\sim \lVert \nabla^2\eta\rVert_{L^2_x(\R^2)}$ and Young's inequality we arrive at
\begin{align*}
\frac12\frac{d}{dt}\int_{\R^2} |\nabla\eta|^2\dx+\frac{\kappa_\ast}{2}\int_{\R^2}|\Delta\eta|^2\dx
\leq  C(\kappa_\ast)\|u\|_{ L^4_x}^4\lVert \nabla\eta\rVert_{ L^2_x}^2.
\end{align*}
Gronwall's inequality    gives
\begin{align*}
 \|\nabla\eta(T)\|_{L^2_x}^2+\| \nabla^2\eta\|_{L^2_T L^2_x}^2
&\leq  C(\kappa_\ast)\|\nabla\eta_0\|_{L^2_x}^2\exp(C(\kappa_\ast)\|u\|_{L^4_T L^4_x}^4)
\end{align*}
for any positive time $T>0$.
Thus by the $\eta$-equation
\begin{align*}
\|\d_t\eta\|_{L^2_T L^2_x} 
&= \|u\cdot\nabla\eta- \kappa\Delta\eta\|_{L^2_T L^2_x}
\leq  \|u\|_{L^4_T L^4_x} \| \nabla\eta \|_{L^4_T L^4_x}+\kappa^\ast \| \Delta\eta\|_{L^2_T L^2_x}
\\
&\leq C(\kappa_\ast,  \kappa^\ast)\|\nabla\eta_0\|_{L^2_x}\exp(C(\kappa_\ast) \|u\|_{L^4_T L^4_x}^4).
\end{align*}

By virtue of the equivalence relation \eqref{eq:eta,theta}:
$$
\|\nabla\theta\|^2_{L^\infty_TL^2_x}+\|\nabla^2 \theta\|^2_{L^2_T L^2_x}
\leq   C(\kappa_\ast,\|a\|_\Lip) (\|\nabla\eta\|_{L^\infty_T L^2_x}^2+ (1+ \|\nabla\eta\|_{L^\infty_T L^2_x}^2)\| \nabla^2\eta\|_{L^2_T L^2_x}^2)  
$$
 and \eqref{theta:LinftyL2}-\eqref{u:L4}, we have the a priori $H^1$-Estimate \eqref{energy:theta:1}  for $\theta$:
\begin{equation}\label{theta:LinftyH1}\begin{split}
&\|\theta\|^2_{L^\infty_TH^1_x}+\|\nabla \theta\|^2_{L^2_T H^1_x} +\|\d_t\theta\|^2_{L^2_TL^2_x}
\\
&\le C(\kappa_\ast,\|a\|_\Lip, \kappa^\ast) \|\theta_0\|^2_{H^1}(1+\|\nabla\theta_0\|^2_{L^2})\exp(C(\kappa_\ast) 
 \|u\|_{L^4_T L^4_x}^4).
\end{split}\end{equation} 

Therefore both the parabolic equations \eqref{theta:eq} and \eqref{eta:eq} for $\theta$ and $\eta$ hold in $L^2_\loc([0,\infty); L^2(\R^2))$.
A standard density argument ensures the $H^1$-Estimate \eqref{energy:theta:1}  for $\theta$ and $\theta\in C([0,\infty); H^1_x(\R^2))$.
\end{paragraph}


\begin{paragraph}{Uniqueness} 
Let  $(\theta_1,u_1,\Pi_1)$ and $(\theta_2,u_2,\Pi_2)$ be two weak solutions of the Cauchy problem \eqref{Boussinesq} with the same initial data   $(\theta_0,u_0) \in H^1(\R^2)\times (L^2(\R^2))^2$, which satisfy the energy estimates \eqref{energy:u:1}-\eqref{energy:theta:1}.
Recall \eqref{eta} for the definition of the function $A$, and we set 
$$\eta_1=A(\theta_1),\quad \eta_2=A(\theta_2).$$
We consider the difference 
$$(\dot \eta,\dot u,\nabla\dot\Pi)=(\eta_1-\eta_2,u_1-u_2,\nabla\Pi_1-\nabla\Pi_2),$$
which lies in 
 \begin{align*}
 &\bigl( C([0,\infty);H^1(\R^2))\cap L^2_\loc([0,\infty);H^2(\R^2))\bigr)\\
&\times\bigl(C([0,\infty); (L^2(\R^2))^2)\cap L^2_\loc([0,\infty); (H^1(\R^2))^2)\bigr)\times L^2_\loc([0,\infty); (H^{-1}(\R^2))^2)\bigr).
\end{align*}
It satisfies the following Cauchy problem
\begin{equation}\label{diff}
\left\{
\begin{aligned}
&\d_t\dot \eta+ u_1\cdot \nabla\dot \eta-\kappa_1\Delta \dot \eta =\dot\kappa\Delta \eta_2 -\dot u\cdot\nabla \eta_2,\\
&\d_t \dot u+u_1\cdot \nabla\dot u-\di(\mu_1 S\dot u)+\nabla\dot \Pi= \dot\theta \vec{e_2}-\dot u\cdot \nabla u_2+\di(\dot\mu Su_2),\\
&\di\dot u=0, \\
&(\dot\eta_0,\dot u_0)=(0,0),
\end{aligned}
\right.
\end{equation}
where 
$$
\kappa_1=a(\theta_1),\,\mu_1=b(\theta_1),\,\dot\theta=\theta_1-\theta_2,\,\dot\kappa=a(\theta_1)-a(\theta_2),\,\dot\mu=b(\theta_1)-b(\theta_2).
$$
Similarly as in \eqref{eq:eta,theta} we have the following equivalence relationships
\begin{equation}\label{eq:dot,eta,theta}\begin{split}
&\kappa_\ast\|\dot\theta\|_{L^2_x}\leq \|\dot\eta\|_{L^2_x}\leq\kappa^\ast\|\dot\theta\|_{L^2_x},\\
& \|\nabla \dot\eta\|_{L^2_x}\leq \|a\|_\Lip\|\nabla\theta_1\|_{L^4_x}\|\dot\theta\|_{L^4_x}
+\kappa^\ast\|\nabla\dot\theta\|_{L^2_x},
\\
&\|\nabla\dot\theta\|_{L^2_x}\leq \frac{\|a\|_\Lip}{\kappa_\ast^3}\|\nabla\eta_1\|_{L^4_x}\|\dot\eta\|_{L^4_x}
+\frac{1}{\kappa_\ast}\|\nabla \dot\eta\|_{L^2_x},
\end{split}\end{equation}   
and correspondingly we have
\begin{equation}\label{eq:mu,eta}
\begin{split}
\|(\dot\kappa, \dot \mu)\|_{H^1_x}
&\leq C(\|(a,b)\|_{\Lip},   \|(a', b')\|_\Lip,\kappa_\ast)(\|\nabla\eta_1\|_{L^4_x}\|\dot\eta\|_{L^4_x}+\|\dot\eta\|_{H^1_x})
\\
&\leq C(\|(a,b)\|_{\Lip},   \|(a', b')\|_\Lip,\kappa_\ast)(1+\|\nabla\eta_1\|_{L^4_x}) \|\dot\eta\|_{H^1_x}.
\end{split}\end{equation}

We are going to sketch the derivation of the $H^1\times L^2$-Estimate for $(\dot\eta, \dot u)$.
\begin{enumerate}[(i)]
\item {$L^2$ estimate of $\dot \eta$.}

We take the $L^2(\R^2)$-inner product between $\eqref{diff}_1$ and $\dot \eta$ to derive
\begin{equation}\label{da}
\begin{aligned}
   &\frac12\frac{d}{dt}\int_{\R^2}|\dot \eta|^2+\int_{\R^2}\kappa_1 |\nabla \dot{\eta}|^2
   \le
   \int_{\R^2}|\dot \eta \nabla\kappa_1\cdot \nabla \dot \eta|+ |\dot \eta\dot u\cdot\nabla \eta_2|+ |\dot \kappa\Delta \eta_2\dot \eta| .
   \end{aligned}
   \end{equation}
   The righthand side can be bounded by 
   \begin{align*}
  & \|\dot \eta\|_{L^4_x} \|\nabla\kappa_1\|_{L^4_x}\|\nabla\dot \eta\|_{L^2_x}
   + \| \nabla \eta_2\|_{L^2_x}\|\dot u\|_{L^4_x}\|\dot \eta\|_{L^4_x}
   +\| \Delta \eta_2\|_{L^2_x}\|\dot \kappa\|_{L^4_x}\|\dot \eta\|_{L^4_x}
   \\
   &\leq C(\|a\|_\Lip)
   \Bigl( \|\nabla\theta_1\|_{L^4_x}\|\dot \eta\|_{L^2_x}^{\frac12} \|\nabla\dot \eta\|_{L^2_x}^{\frac32}
   + \| \nabla \eta_2\|_{L^2_x}\|\dot u\|_{L^2_x}^{\frac12}\|\nabla\dot u\|_{L^2_x}^{\frac12}\|\dot \eta\|_{L^2_x}^{\frac12}\|\nabla\dot \eta\|_{L^2_x}^{\frac12}
   \\
   &\qquad 
  + \| \Delta \eta_2\|_{L^2_x}\|\dot \theta\|_{L^2_x}^{\frac12}\|\nabla\dot \theta\|_{L^2_x}^{\frac12}\|\dot \eta\|_{L^2_x}^{\frac12}\|\nabla\dot \eta\|_{L^2_x}^{\frac12}\Bigr)
  \\
  &\leq \frac{\kappa_\ast}{2}\|\nabla\dot \eta\|_{L^2_x}^2+\frac{\mu_\ast}{4}\|\nabla\dot u\|_{L^2_x}^2
  \\
  &\quad+ C(\|a\|_\Lip, \kappa_\ast,  \mu_\ast)
  \Bigl(\|\nabla\theta_1\|_{L^4_x}^4 + \| \nabla \eta_2\|_{L^2_x}^2+  \| \Delta \eta_2\|_{L^2_x}^2\Bigr)
  \times\bigl(\|\dot \eta\|_{L^2_x}^{2}+\|\dot u\|_{L^2_x}^{2} \bigr).
   \end{align*}     

\item {$L^2$ estimate of $\dot u$.} 

We take the $L^2$ inner product of the equation $\eqref{diff}_2$ and $\dot u$ to derive
\begin{equation}\label{dau}
\begin{aligned}
  &\frac12\frac{d}{dt}\int_{\R^2}|\dot u|^2+\frac12\int_{\R^2}\mu_1|S \dot u|^2
  \le \int_{\R^2} |\dot\theta \dot u|+ |\dot u|^2|\nabla u_2|+ | \dot\mu S u_2:\nabla \dot u|.
  \end{aligned}
  \end{equation}
  The righthand side can be bounded by
  \begin{align*}
  \|\dot u\|_{L^2_x}\|\dot \theta\|_{L^2_x}
  +\|\dot u\|_{L^4_x}^2\|\nabla u_2\|_{L^2_x}
  +\|\dot \mu\|_{H^1_x}\|S u_2:\nabla\dot u\|_{H^{-1}_x}
  \end{align*}
  which, by use of the Sobolev embedding  $L^1(\R^2)\hookrightarrow H^{-1}(\R^2) $, is bounded by
  \begin{align*}
  &\frac{\mu_\ast}{4}\|\nabla\dot u\|_{L^2_x}^2 
  \\
&  +C(\|a\|_\Lip,\kappa_\ast, \mu_\ast)(1+\|\nabla u_2\|_{L^2_x}^2) (\|\dot u\|_{L^2_x}^2+\|\dot \eta\|_{L^2_x}^2)
  +C(\mu_\ast)\|\nabla u_2\|_{L^2_x}^2\|\dot\mu\|_{H^1_x}^2.  
  \end{align*} 
 

\item {$L^2$ estimate of $\nabla \dot \eta$.}

We take the $L^2$ inner product of the equation $\eqref{diff}_1$ and $\Delta\dot \eta$ to derive
\begin{align*}
    &\frac{1}{2}\frac{d}{dt}\int_{\R^2}|\nabla\dot \eta|^2+\int_{\R^2}\kappa_1 |\Delta\dot \eta|^2
    \le \int_{\R^2}| u_1\cdot \nabla\dot \eta\Delta\dot \eta|
    + |\dot u\cdot\nabla \eta_2\Delta\dot \eta|+ |\dot\kappa\Delta \eta_2\Delta\dot \eta| 
\end{align*}
By $L^1(\R^2)\hookrightarrow H^{-1}(\R^2) $ again, the righthand side is bounded similarly by
\begin{align*}
&\frac{\kappa_\ast}{2}\|\Delta\dot\eta\|_{L^2_x}^2+\frac{\mu_\ast}{4}\|\nabla\dot u\|_{L^2_x}^2 
+
C(\kappa_\ast, \mu_\ast) ( \|u_1\|_{L^4_x}^4+\|\nabla \eta_2\|_{L^4_x}^4) 
(\|\dot u\|_{L^2_x}^2+\|\nabla\dot \eta\|_{L^2_x}^2)
\\
&+C(\kappa_\ast)\|\Delta \eta_2\|_{L^2_x}^2\|\dot\kappa\|_{H^{1}_x}^2.
\end{align*} 
\end{enumerate}

To conclude, by virtue of the above estimates and \eqref{eq:mu,eta}, we have the following inequlity
\begin{equation*}
\begin{aligned}
& \frac{d}{dt}(\|\dot \eta \|_{H^1_x}^2+\|\dot u\|_{L^2_x}^2)+\|\nabla \dot u\|_{L^2_x}^2+\|\nabla\dot \eta\|_{H^1_x}^2  
\\
&\le C(\|(a,b)\|_{\Lip},   \|(a', b')\|_\Lip, \kappa_\ast, \mu_\ast) B(t)  (\|\dot u\|_{L^2_x}^2+\|\dot \eta\|_{H^1_x}^2),
\end{aligned}
\end{equation*}
where 
\begin{align*}
B(t)&=\Bigl(\|\nabla\theta_1\|_{L^4_x}^4 + \| \nabla \eta_2\|_{L^2_x}^2+  \| \Delta \eta_2\|_{L^2_x}^2+1+\|\nabla u_2\|_{L^2_x}^2 +\|u_1\|_{L^4_x}^4+\|\nabla \eta_2\|_{L^4_x}^4\Bigr)
\\
&\quad\times (1+\|\nabla\eta_1\|_{L^4_x})\in L^1_\loc([0,\infty)).
\end{align*}
Gronwall's inequality implies then $\dot\eta=0$ and $\dot u=0$.
The uniqueness of the weak solutions follows. 
\end{paragraph}


\subsection{Propagation of the general $H^s$-regularities}\label{sec:fractional}
After the derivation of the a priori $H^s_x$, $s\in (0,2)$-Estimate for a general linear parabolic equation in Subsection \ref{subs:parabolic}, we are going to derive the precise $H^s_x$-estimates \eqref{theta:Hsnu}-\eqref{u:Hs3} in Remark \ref{rem:Hs} in the subsequent subsections:
\begin{itemize}
\item In Subsection \ref{subs:horizontal} the global-in-time $H^s_x(\R^2)\times (L^2_x(\R^2))^2$, $s\in (1,2]$ -regularities (i.e. \eqref{theta:Hsnu}) will be established, where the endpoint case $s=2$ (i.e. \eqref{theta:Hs2}) will be treated separately. 
\item In Subsection \ref{subs:vertical} the global-in-time $H^1_x(\R^2)\times (H^s_x(\R^2))^2$, $s\in (0,2]$-regularities (i.e. \eqref{u:Hsnu}) will be established, where the endpoint case $s=2$ (i.e. \eqref{u:Hs2}) will be treated separately.
\item In Subsection \ref{subs:s:2} the global-in-time  $  H^{s}_x(\R^2)\times (H^{s-2}_x(\R^2))^2$ (i.e. \eqref{theta:Hs23}-\eqref{theta:Hs3}) and 
  $ H^{s-1}_x\times (H^s_x(\R^2))^2$, $s>2$-regularities (i.e. \eqref{u:Hs23}-\eqref{u:Hs3}) will be established respectively. 
\end{itemize}
As far as the borderline estimates  \eqref{theta:Hsnu}-\eqref{u:Hs3} are established, the global-in-time regularity \eqref{regularity} follows immediately.

For readers' convenience we recall here briefly the Littlewood-Paley dyadic decomposition and the $H^s(\R^n)$-norms (see e.g.  Chapter 2 in the book \cite{BCD} for more details).
We fix a nonincreasing radial function $\chi\in C^\infty_c(B_{\frac{4}{3}})$ with $\chi(x)=1$ with $x\in B_1$, where $B_r\subset\R^n$ denotes the ball centered at $0$ with radius $r$.   We define the function $\varphi(\xi)=\chi(\frac{\xi}{2})-\chi(\xi)$ and $\varphi_j(\xi)=\varphi(2^{-j}\xi)$ with $j\ge 0$. We do the Littlewood-Paley decomposition in the following way  
\begin{align}\label{LP}
   g=\Delta_{-1}g+\sum_{j\ge 0} \Delta_j g,
\end{align}
where 
\begin{align*}
    \cF(\Delta_{-1}g)(\xi)=\chi(\xi)\cF(g)(\xi),\quad \cF(\Delta_{j}g)(\xi)=\varphi_j(\xi)\cF(g)(\xi), \quad j\geq0,
\end{align*}
and $\cF$ denotes the Fourier transform. 
We have the following Bernstein's inequalities for some universal constant $C$ (depending only on $n$)
\begin{align}\label{Bernstein}
&\|\Delta_{-1}g\|_{L^2(\R^n)}\le C\|g\|_{L^2(\R^n)},   \notag
\\
&C^{-1}2^{j}\|\Delta_j g\|_{L^2(\R^n)}\le \|\nabla (\Delta_j g)\|_{L^2(\R^n)}\leq C2^{j}\|\Delta_j g\|_{L^2(\R^n)},\quad 
\forall j\geq0.
\end{align} 
Let $s\ge0$ and $p,\,r\ge1$. We define the nonhomogeneous Besov spaces $B^s_{p,r}(\R^n)$ as the spaces consisting of all tempered distributions $g\in\cS'(\R^n)$ satisfying
\begin{equation*}
    \|g\|_{B^s_{p,r}(\R^n)}=\big\|(2 ^{js}\|\Delta_j g\|_{L^p(\R^n)})_{j\ge-1}\big\|_{l^r}<\infty.
\end{equation*}
The inhomogeneous Sobolev spaces $H^s(\R^n)=   B^s_{2,2}(\R^n)$ can be defined by 
\begin{equation*}
H^s(\R^n)=\{g\in\cS'(\R^n)\mid 
\|g\|_{H^s(\R^n)}=\Bigl( \int_{\R^n}(1+|\xi|^2)^{\frac{s}{2}}|\cF(g)(\xi)|^2\,d\xi\Bigr)^{1/2}<\infty\}, 
\end{equation*}
where the $H^s(\R^n)$-norm reads in terms of Littlewood-Paley decomposition as follows
\begin{equation}\label{Hs:LP} 
\|g\|_{H^s(\R^n)}\sim \|g\|_{L^2(\R^n)}+\Bigl(\sum_{j\geq 0}2^{2js}\|\Delta_j g\|_{L^2(\R^n)}^2\Bigr)^{\frac12}.
\end{equation} 
It is straightforward to derive the following interpolation inequality  
   \begin{equation}\label{interpolation}
      \|u\|_{H^{t_\sigma}}\le C\|u\|_{H^{t_0}}^{1-\sigma}\|u\|_{H^{t_1}}^\sigma,
      \hbox{ where } t_\sigma=(1-\sigma)t_0+\sigma t_1 ,  \sigma\in[0,1].
   \end{equation}
   
   We are going to use the following known  estimates  to control the nonlinear terms in the Boussinesq system \eqref{BS}.
  \begin{lemma} 
    \label{Liao}
We have the following commutator, product and composition estimates.
\begin{enumerate}[(1)]
\item {  \cite[Proposition 2.4]{Danchin:Liao}}   In the low regularity regime where  $(s,\nu)\in\R^2$ satisfy
\begin{equation*}
    -1<s<\nu+1,\quad\text{and } -1<\nu<1, 
\end{equation*}
the following commutator estimate holds true for some constant $C$ depending on $s,\nu$ and for the spacial dimension two
\begin{equation}\label{commutator:1}
    \|\big(2^{js}\|[\phi,\Delta_j]\nabla\psi\|_{L^{2}(\R^2)}\big)_{j\ge 1}\|_{l^1}\le C\|\nabla\phi\|_{H^{\nu}(\R^2)}\|\nabla\psi\|_{H^{s-\nu}(\R^2)}.
\end{equation}
\item{\cite[Lemma 2.100]{BCD}} For any $s>0$, the following commutator estimate holds true
\begin{equation}\label{commutator:2}\begin{split}
    &\|\big(2^{js}\|[\phi,\Delta_j]\nabla\psi\|_{L^{2}(\R^n)}\big)_{j\ge 1}\|_{l^2}
    \\
&    \le C\bigl( \|\nabla\phi\|_{L^\infty(\R^n)}\|\nabla\psi\|_{H^{s-1}(\R^n)}+\|\nabla\phi\|_{H^{s-1}(\R^n)}\|\nabla\psi\|_{L^\infty(\R^n)}).
\end{split}\end{equation}

\item {\cite[Corollary 2.86]{BCD}} For any $s>0$, the following product estimate holds true
\begin{equation}\label{product}
\|\phi\psi\|_{H^{s}(\R^n)}\leq  C\bigl( \|\phi\|_{L^\infty(\R^n)}\|\psi\|_{H^{s}(\R^n)}+\|\phi\|_{H^{s}(\R^n)}\|\psi\|_{L^\infty(\R^n)}).
\end{equation}

\item {\cite[Theorem 2.87 \& Theorem 2.89]{BCD}} For any $s>0$ and  $g\in C^{k+1}$ with $k=[s]\in \N$, the following composition estimate holds true  
\begin{equation}\label{composition0}
\|\nabla (g\circ \theta)\|_{H^{s-1}(\R^n)}\leq C(g, \|\theta\|_{L^\infty(\R^n)})\|\nabla\theta\|_{H^{s-1}(\R^n)}.
\end{equation}
If $g\in C^{k+1}_b$ with $k=[s]\in \N$, then the above estimate can be improved in the spacial dimension two as follows
\begin{equation}\label{composition}
\|\nabla (g\circ \theta)\|_{H^{s-1}(\R^2)}\leq C(\|g\|_{C^{k+1}}, \|\theta\|_{H^1(\R^2)})\|\nabla\theta\|_{H^{s-1}(\R^2)}.
\end{equation}

\end{enumerate}
\end{lemma}
The commutator estimate \eqref{commutator:1} will present its power in the low regularity regime (see Subsections \ref{subs:parabolic}-\ref{subs:vertical} below),  and the classical commutator estimate \eqref{commutator:2} will help in the high regularity regime (see Subsection \ref{subs:s:2}).

The composition estimate \eqref{composition} will help to bound the diffusion coefficients $\kappa, \mu$ in terms of $\theta$ in the low regularity regime, where only  $H^1(\R^2)$-norm (instead of $L^\infty_x$-norm) of $\theta$ is available.

\subsubsection{Estimates for the general parabolic equations}\label{subs:parabolic} 
We derive in this paragraph a priori $H^s$, $s\in (0,2)$-Estimates for a general   linear parabolic equation, which should be of independent interest. 
\begin{lemma}\label{lem:parabolic}
Let $\psi=\psi(t,x): [0,\infty)\times\R^2\mapsto \R^m$, $m\geq 1$ be a smooth solution with sufficiently decay of the following linear parabolic equation
\begin{equation}\label{parabolicc}
\left\{
\begin{aligned}
&\d_t\psi+ u\cdot \nabla_x \psi-\di_x(\kappa \nabla_x \psi)=f,\\
&\psi|_{t=0}=\psi_0,
\end{aligned}
\right.
\end{equation}
where 
\begin{itemize}
\item $u=u(t,x): \R^+\times\R^2\mapsto \R^2$ is a given divergence-free vector field: $\div_x u=0$;
\item $\kappa=\kappa(t,x):\R^+\times\R^2\to[\kappa_*,\kappa^*]$ with $\kappa_*,\kappa^*\in(0,\infty)$;
\item  $f=f(t,x): \R^+\times\R^2\mapsto \R^m$ denotes the given external force. 
\end{itemize}

Then the following a priori $H^s_x$-Estimates  for \eqref{parabolicc} holds true:
\begin{equation}\label{psi:Hs}\begin{split}
&\|\psi\|_{L^\infty_T H^s_x}^2+\|\nabla\psi\|_{L^2_T H^{s}_x}^2
\leq C(\kappa_\ast)\Bigl( \|\psi_0\|_{H^s_x}^2+\|f\|_{L^2_T H^{s-1}_x}^2  \Bigr)\times
\\
&\times \exp\Bigl(C(\kappa_\ast,s,\nu)(\|\nabla u\|_{L^2_T L^2_x}^2+\|\nabla\kappa\|_{L^{\frac2\nu}_T H^\nu_x}^{2/\nu}+ \|f\|_{L^1_T H^{-s}_x})\Bigr) \\
&\qquad\qquad\qquad
\hbox{ for any }s\in(0,2) \hbox{ with }\nu\in (s-1,1).
\end{split}\end{equation} 
\end{lemma}
\begin{proof} We are going to derive the 
It is straightforward to derive the following  $L^2_x$-Estimate  by simply taking the $L^2(\R^2)$ inner product of the equation \eqref{parabolicc} and $\psi$  itself
 \begin{equation}\label{L2:psi}
\|\psi\|_{L^\infty_TL^2_x}^2+\|\nabla \psi\|^2_{L^2_TL^2_x}\le 
C(\kappa_\ast) \Bigl( \|\psi_0\|_{L^2_x}^2+\int^T_0 \langle \psi, f\rangle_{H^{s}_x, H^{-s}_x}\,dt\Bigr),
\quad\forall s\in\R.
\end{equation} 

We next consider the a priori estimates for  the $H^s(\R^2)$-norm. 
By virtue of the  description \eqref{Hs:LP} of the $H^s(\R^2)$-norm,
 we consider the dyadic piece of $\psi$:
\begin{equation}\label{Deltaj:psi}
\psi_j:=\Delta_j\psi, \quad j\geq 0.
\end{equation}
where the operator $\Delta_j$ is defined in \eqref{LP}.
We apply $\Delta_j$ to the linear $\psi$-equation to derive the equation for $\psi_j$:
    \begin{equation}\label{theta:j:s}
        \d_t\psi_j+u\cdot\nabla\psi_j-\di(\kappa\nabla\psi_j)=[u,\Delta_j]\cdot\nabla\psi-\div([\kappa,\Delta_j]\nabla\psi)+f_j,
        \quad j\geq 0.
    \end{equation} 
 We take the $L^2$ inner product of the equation \eqref{theta:j:s} and $\psi_j$ and make use of $\div u=0$ and $\kappa\geq\kappa_\ast$ to derive
    \begin{equation*}\label{theta:j:l2:1}
    \begin{aligned}
    &\frac12\frac{d}{dt}\|\psi_j\|_{L^2_x}^2+\kappa_*\|\nabla\psi_j\|_{L^2_x}^2 
    \le   \|\psi_j\|_{L^2_x}\|[u,\Delta_j]\cdot\nabla\psi\|_{L^2_x}
    \\
    &\qquad\qquad+\|\nabla\psi_j\|_{L^2_x}\|[\kappa,\Delta_j]\nabla\psi\|_{L^2_x}+\|f_j\|_{L^2_x}\|\psi_j\|_{L^2_x},\, j\geq 0.
    \end{aligned}
    \end{equation*}
  By use of Bernstein's inequality \eqref{Bernstein} we have
    \begin{equation*}\label{theta:j:l2:2}
    \begin{aligned}
     &\frac{d}{dt}\|\psi_j\|_{L^2_x}^2+2^{2j}\|\psi_j\|_{L^2_x}^2\\
     \le & C(\kappa_\ast)\|\psi_j\|_{L^2_x}\left(\|[u,\Delta_j]\cdot\nabla\psi\|_{L^2_x}+2^j\|[\kappa,\Delta_j]\nabla\psi\|_{L^2_x}+\|f_j\|_{L^2_x}\right),
    \end{aligned}
    \end{equation*}
    that is, 
    \begin{equation}\label{theta:j:l2:3}
    \begin{aligned}
     &\frac{d}{dt}\|\psi_j\|_{L^2_x}+2^{2j}\|\psi_j\|_{L^2_x}\\
     \le & C(\kappa_\ast) \left(\|[u,\Delta_j]\cdot\nabla\psi\|_{L^2_x}+2^j\|[\kappa,\Delta_j]\nabla\psi\|_{L^2_x}+\|f_j\|_{L^2_x}\right),
     \quad j\geq 0.
    \end{aligned}
    \end{equation}
    
    We   make use of the commutator estimate \eqref{commutator:1} in Lemma \ref{Liao} to estimate the commutators $\|[u,\Delta_j]\cdot\nabla\psi\|_{L^2_x}$ and $2^j\|[\kappa,\Delta_j]\nabla\psi_j\|_{L^2_x}$ in the above inequality in the following way.
     Let $(l_j)_{j\ge 0}$ be a normalised sequence in $\ell^1(\N)$ such that $l_j\geq 0$ and $\sum_{j\ge 0}l_j=1$.  
     Then we have 
     \begin{equation}\label{comm:est}\begin{split}
     &\|[u,\Delta_j]\nabla\psi\|_{L^{2}} \le C(s) l_j 2^{j(1-s)}\|\nabla u\|_{L^2_x}\|\nabla\psi\|_{H^{s-1}_x},
     \hbox{ for } s\in (0,2),\\
 &2^j\|[\kappa,\Delta_j]\nabla\psi\|_{L^2_x}
     \le C(s,\nu) l_j 2^{j(1-s)}\|\nabla \kappa\|_{H^{\nu}_x}\|\nabla\psi\|_{H^{s-\nu}_x}
     \\
&\qquad\qquad\qquad\qquad\qquad
     \hbox{ for } \nu\in (-1,1),\, s\in (-1,\nu+1).
     \end{split}\end{equation}  
     Therefore we have
    \begin{align*}
       &\frac{d}{dt}\|\psi_j\|_{L^2_x}+2^{2j}\|\psi_j\|_{L^2_x}
       \\
   &    \leq 
  C(\kappa_\ast,s,\nu) l_j 2^{j(1-s)}\left(\|\nabla u\|_{L^2_x}\|\nabla\psi\|_{H^{s-1}_x}+\|\nabla \kappa\|_{H^{\nu}_x}\|\nabla\psi\|_{H^{s-\nu}_x}\right)+C(\kappa_\ast)\|f_j\|_{L^2_x}\\ 
       &  \qquad\qquad\qquad\qquad\qquad\qquad 
       \hbox{ for }\nu\in (-1,1),\,s\in (0,\nu+1),\,j\geq 0.
    \end{align*}
We use Duhamel's Principle to derive 
\begin{equation}\label{theta:k:hs:0}
\begin{aligned}
  &  \|\psi_j\|_{L^2_x} \le e^{-  t2^{2j}}\|(\psi_0)_j\|_{L^2_x}+C(\kappa_\ast) \int_0^t e^{- (t-\tau)2^{2j}}\|f_j(\tau)\|_{L^2_x}\,d\tau \\
    &\quad+C(\kappa_\ast,s,\nu) 2^{j(1-s)}l_j\int_0^t e^{- (t-\tau)2^{2j}} \big(\|\nabla u(\tau)\|_{L^2_x}\|\nabla\psi(\tau)\|_{H^{s-1}_x}\\
    &\qquad\qquad\qquad
    +\|\nabla\kappa(\tau)\|_{H^{\nu}_x}\|\nabla\psi(\tau)\|_{H^{s-\nu}_x}\big)\,d\tau,\quad j\ge 0. 
    \end{aligned}
\end{equation}

We multiply the inequality \eqref{theta:k:hs:0} by $2^{js}$ to derive
\begin{equation}\label{psi:j:s}
\begin{aligned}
   & 2^{js}\|\psi_j\|_{L^2_x} \le 2^{js}e^{-  t2^{2j}}\|(\psi_0)_j\|_{L^2_x}+C(\kappa_\ast) 2^{js}\int_0^t e^{- (t-\tau)2^{2j}} \|f_j\|_{L^2_x}\,d\tau\\
    &\quad+C(\kappa_\ast,s,\nu) 2^{j}l_j\int_0^t e^{-(t-\tau)2^{2j}} \big(\|\nabla u(\tau)\|_{L^2_x}\|\nabla\psi(\tau)\|_{H^{s-1}_x}\\
    &\qquad\qquad\qquad
    +\|\nabla\kappa(\tau)\|_{H^{\nu}_x}\|\nabla\psi(\tau)\|_{H^{s-\nu}_x}\big)\,d\tau,\quad j\ge 0. 
\end{aligned}
\end{equation}
We take $L^\infty([0,T])$-norm in $t$ of \eqref{psi:j:s} and the $L^2([0,T])$-norm in $t$ of $2^j\cdot$\eqref{psi:j:s}, to derive by use of Young's inequality that
\begin{equation}\label{psi:j:s:t}
\begin{aligned}
   & 2^{js}\|\psi_j\|_{L^\infty_T L^2_x}+2^{j(s+1)}\|\psi_j\|_{L^2_T L^2_x} \le 2^{js}\|(\psi_0)_j\|_{L^2_x}+C(\kappa_\ast) 2^{j(s-1)}\|f_j\|_{L^2_T L^2_x}\\
    &\quad+C(\kappa_\ast,s,\nu) l_j \Bigl\|  \|\nabla u\|_{L^2_x}\|\nabla\psi\|_{H^{s-1}_x}+\|\nabla \kappa\|_{H^{\nu}_x}\|\nabla\psi\|_{H^{s-\nu}_x}\Bigr\|_{L^2_T}.
\end{aligned}    
\end{equation}
We take square of \eqref{psi:j:s:t} and sum them up for $j\in\N$ to derive
\begin{align*}
&\sum_{j\geq 0}\Bigl( 2^{2js}\|\psi_j\|_{L^\infty_T L^2_x}^2+2^{2j(s+1)}\|\psi_j\|_{L^2_T L^2_x}^2\Bigr)\\
&\lesssim_{\kappa_\ast,s,\nu} \sum_{j\geq 0}\Bigl(2^{2js}\|(\psi_0)_j\|_{L^2_x}^2+  2^{2j(s-1)}\|f_j\|_{L^2_T L^2_x}^2\Bigr)\\
&+\int^T_0 \|\nabla u\|_{L^2_x}^2\|\nabla\psi\|_{H^{s-1}_x}^2+\|\nabla\kappa\|_{H^\nu_x}^2\|\nabla\psi\|_{H^{s-\nu}_x}^2 dt,\quad j\geq 0,
\end{align*}
that is, by virtue of the $L^2$-estimate \eqref{L2:psi},
\begin{align*}
&\|\psi\|_{L^\infty_T H^s_x}^2+\|\nabla\psi\|_{L^2_T H^{s}_x}^2
\lesssim_{\kappa_\ast,s,\nu}  \Bigl( \|\psi_0\|_{H^s_x}^2+\|f\|_{L^2_T H^{s-1}_x}^2+\int^T_0 \|\psi\|_{H^s_x}\|f\|_{H^{-s}_x} dt\\
&\qquad 
+\int^T_0 \|\nabla u\|_{L^2_x}^2\|\nabla\psi\|_{H^{s-1}_x}^2+\|\nabla\kappa\|_{H^\nu_x}^2\|\nabla\psi\|_{H^{s-\nu}_x}^2 dt\Bigr).
\end{align*}
We next consider the norm $\|\nabla\psi\|_{H^{s-\nu}_x}$.
By the interpolation inequality  \eqref{interpolation} we have
$$
\|\nabla\psi\|_{H^{s-\nu}_x} \leq C\|\nabla \psi\|_{H^{s-1}_x}^{\nu}\|\nabla \psi\|_{H^s_x}^{1-\nu},
\quad \nu\in (0,1),
$$
which implies by Young's inequality that
\begin{align*}
\int^T_0 \|\nabla\kappa\|_{H^\nu_x}^2\|\nabla\psi\|_{H^{s-\nu}_x}^2 dt
&\leq \int^T_0 \|\nabla\kappa\|_{H^\nu_x}^2\|\nabla \psi\|_{H^{s-1}_x}^{2\nu}\|\nabla \psi\|_{H^s_x}^{2(1-\nu)} dt
\\
&\leq \varepsilon \|\nabla\psi\|_{L^2_T H^s_x}^2+C_{\varepsilon}\int^T_0 \|\nabla\kappa\|_{H^\nu_x}^{2/\nu}\|\nabla \psi\|_{H^{s-1}_x}^{2}\,dt.
\end{align*}
To conclude,  by taking $\varepsilon$ small enough and  Gronwall's inequality, we derive the $H^s$-Estimate \eqref{psi:Hs}.
\end{proof}
 

\subsubsection{\texorpdfstring{Case $(\theta_0,u_0)\in H^{s}(\R^2)\times (L^2(\R^2))^2$, $s\in (1,2]$}{} } \label{subs:horizontal}
In this subsection we are going to prove the $H^s$-Estimates \eqref{theta:Hsnu} for  the unique  solution  $(\theta, u)$ of the Boussinesq equations \eqref{BS} with the initial data $(\theta_0, u_0)\in H^{s}(\R^2)\times (L^2(\R^2))^2$, $s\in (1,2)$, following exactly the procedure   in Subsection \ref{subs:parabolic}.  
We will pay more attention on the ``nonlinearities'' in the equations such as $\kappa=a(\theta)$, $u\cdot\nabla u$ when using the commutator estimates and will sketch the proof.

The endpoint estimate \eqref{theta:Hs2} for $(\theta_0, u_0)\in H^{2}(\R^2)\times (L^2(\R^2))^2$ follows similarly as in the proof for the $H^1$-Estimate for $\theta$ in Subsection \ref{subs:uniqueness} and we will sketch the proof.

\begin{paragraph}{Case $(\theta_0,u_0)\in H^{s}(\R^2)\times (L^2(\R^2))^2$, $s\in (1,2)$} 
Similarly as   \eqref{theta:j:l2:3}, we have the following preliminary estimate for $\theta_j=\Delta_j \theta$:
    \begin{equation}\label{theta:j}
       \frac{d}{dt}\|\theta_j\|_{L^2_x}+2^{2j}\|\theta_j\|_{L^2_x}\le C(\kappa_\ast) \left(\|[u,\Delta_j]\cdot\nabla\theta\|_{L^2_x}+2^j\|[\kappa,\Delta_j]\nabla\theta\|_{L^2_x}\right),\quad j\geq 0.
    \end{equation} 
    By use of the commutator estimates \eqref{comm:est} and the action estimate \eqref{composition}: 
    $$\|\nabla\kappa\|_{H^\nu}\leq C(\|a\|_{C^2}, \|\theta\|_{H^1})\|\nabla\theta\|_{H^\nu}\hbox{ for }\nu\in (0,1),$$
     we derive similar as \eqref{psi:j:s:t}
\begin{equation}\label{theta:j:s:t}
\begin{aligned}
   & 2^{2js}\|\theta_j\|_{L^\infty_T L^2_x}^2+2^{2j(s+1)}\|\theta_j\|_{L^2_T L^2_x}^2 \le 2^{2js}\|(\theta_0)_j\|_{L^2_x}^2 \\
    & +C(\kappa_\ast, s, \nu, \|a\|_{C^2},\|\theta\|_{L^\infty_T H^1_x}) (l_j)^2 \int^T_0 \Bigl( \|\nabla u\|_{L^2_x}^2\|\nabla\theta\|_{H^{s-1}_x}^2
    \\
    &\qquad\qquad
    +\|\nabla \theta\|_{H^{\nu}_x}^2\|\nabla\theta\|_{H^{s-\nu}_x}^2\Bigr) dt,
    \quad 1<s<\nu+1<2.
\end{aligned}    
\end{equation} 
By using the interpolation inequality \eqref{interpolation}, we have
$$
\|\nabla \theta\|_{H^{\nu}_x}\|\nabla\theta\|_{H^{s-\nu}_x}
\leq C\|\nabla\theta\|_{L^2_x}^{1-\nu}\|\nabla\theta\|_{H^1_x}^\nu \|\nabla\theta\|_{H^{s-1}_x}^{\nu}\|\nabla\theta\|_{H^{s}_x}^{1-\nu},
\quad 0< \nu<1.
$$
Recall the  $L^2$-Estimate \eqref{theta:LinftyL2} for $\theta$: 
\begin{equation}\label{theta:L2}
\|\theta\|_{L^\infty_T L^2_x}^2+\|\nabla\theta\|_{L^2_T L^2_x}^2\leq C(\kappa_\ast) \|\theta_0\|_{L^2_x}^2.
\end{equation}
Therefore by Young's inequality we arrive at
\begin{align*}
&\|\theta\|_{L^\infty_T H^s_x}^2+\|\nabla\theta\|_{L^2_T H^{s}_x}^2 \leq   C(\kappa_\ast)\|\theta_0\|_{H^s_x}^2\\
&+C(\kappa_\ast, s, \nu, \|a\|_{C^2},\|\theta\|_{L^\infty_T H^1_x})
\int^T_0 \Bigl(\|\nabla u\|_{L^2_x}^2+ \|\nabla\theta\|_{H^1_x}^2\Bigr)\|\nabla\theta\|_{H^{s-1}_x}^2 dt,
\end{align*}
which, together with Gronwall's inequality, implies \eqref{theta:Hsnu}.
\end{paragraph}

\bigbreak

\begin{paragraph}{Endpoint case $(\theta_0, u_0)\in H^2(\R^2)\times (L^2(\R^2))^2$}
We recall the function  $\eta=A^{-1}(\theta)$ defined in \eqref{eta,theta}, and the parabolic $\eta$-equation \eqref{eta:eq}:
\begin{equation}\label{eta:eq+}
\d_t\eta+u\cdot\nabla\eta-\kappa\Delta\eta=0.
\end{equation}
We are going to derive the a priori $H^2$-Estimate for $\eta$ under the conditions 
\begin{align*}
\div u=0,\quad \nabla u\in L^2_\loc([0,\infty); (L^2(\R^2))^2)\hbox{ and }\nabla\kappa\in L^4_\loc([0,\infty); (L^4(\R^2))^2).
\end{align*}
We test the above $\eta$-equation \eqref{eta:eq+} by $\Delta^2\eta$, to arrive at
\begin{align*}
\frac12\frac{d}{dt}\int_{\R^2}|\Delta\eta|^2\dx
+\int_{\R^2}\kappa|\nabla\Delta\eta|^2\dx
=-\int_{\R^2} \Bigl( u\cdot\nabla\eta \Delta^2\eta+\nabla\kappa\cdot\nabla\Delta\eta \Delta\eta\Bigr)\dx.
\end{align*}
By integration by parts, $\div u=0$ and the embedding $L^1(\R^2)\hookrightarrow H^{-1}(\R^2)$, we derive
\begin{align*}
-\int_{\R^2} u\cdot\nabla\eta \Delta^2\eta \dx
&=\int_{\R^2}\nabla\Delta\eta\cdot\nabla u\cdot\nabla\eta   - \nabla u:\nabla^2\eta \Delta\eta \dx 
\\
&\leq  4\|\nabla u\|_{L^2_x}( \|\nabla\Delta\eta\|_{L^2_x}\|\nabla\eta\|_{H^1_x}+\|\nabla^2\eta\|_{L^4_x}^2)
\\
&\leq \frac{\kappa_\ast}{4}\|\nabla\Delta\eta\|_{L^2_x}^2
+C(\kappa_\ast)\|\nabla u\|_{L^2_x}^2\|\nabla\eta\|_{H^1_x}^2.
\end{align*}
Similarly we have
\begin{align*}
-\int_{\R^2} \nabla\kappa\cdot\nabla\Delta\eta \Delta\eta\dx
&\leq \|\nabla\kappa\|_{L^4_x}\|\nabla\Delta\eta\|_{L^2_x}\|\Delta\eta\|_{L^4_x}
\\
&\leq \frac{\kappa_\ast}{4}\|\nabla\Delta\eta\|_{L^2_x}^2
+C(\kappa_\ast)\|\nabla \kappa\|_{L^4_x}^4\|\Delta\eta\|_{L^2_x}^2.
\end{align*}
To conclude, we have the following a priori $\dot H^2_x$-Estimate for $\eta$ and any positive time $T>0$ by Gronwall's inequality 
\begin{align*}
 \|\Delta\eta(T)\|_{L^2_x}^2+\| \nabla\Delta\eta\|_{L^2_T L^2_x}^2
&\leq  C(\kappa_\ast)(\|\Delta\eta_0\|_{L^2_x}^2+\|\nabla u\|_{L^2_T L^2_x}^2\|\nabla\eta\|_{L^\infty_TL^2_x}^2)
\\
&\qquad\times \exp\Bigl(C(\kappa_\ast)(\|\nabla u\|_{L^2_T L^2_x}^2+\|\nabla\kappa \|_{L^4_T L^4_x}^4)\Bigr).
\end{align*} 
By view of the equivalence relation \eqref{eq:eta,theta} as well as\footnote{It is also straightforward to calculate
\begin{align*}
&\d_{jkl}\eta
=a''(\theta)(\d_j\theta\d_k\theta\d_l\theta)
+a'(\theta)(\d_{jk}\theta\d_l\theta+\d_{jl}\d_k\theta+\d_{kl}\theta\d_j\theta)
+a(\theta)\d_{jkl}\theta,
\\
&\d_{jkl}\theta
=\Bigl(-\frac{a''(A^{-1}(\eta))}{a^4(A^{-1}(\eta))}+\frac{3(a'(A^{-1}(\eta)))^2}{a^5(A^{-1}(\eta))}\Bigr)\d_j\eta\d_k\eta\d_l\eta
\\
&\qquad\quad
-\frac{a'(A^{-1}(\eta))}{a^3(A^{-1}(\eta))}\Bigl(\d_{jk}\eta\d_l\eta+\d_{jl}\eta\d_l\eta+\d_{kl}\eta\d_j\eta\Bigr)
+\frac{1}{a(A^{-1}(\eta))}\d_{jkl}\eta.
\end{align*}}
$$
\|\nabla^3\theta\|_{L^2_TL^2_x}\leq C(\kappa_\ast, \|a\|_{C^2})\Bigl( (\|\nabla\eta\|_{L^4_TL^4_x}^2+\|\nabla^2\eta\|_{L^2_TL^2_x})\|\nabla\eta\|_{L^\infty_TL^2_x}+\|\nabla^3\eta\|_{L^2_TL^2_x}\Bigr),
$$
we derive the $H^2$-Estimate \eqref{theta:Hs2} for $\theta$ by virtue of the $H^1$-Estimate \eqref{theta:LinftyH1}:
\begin{align*}
&\|\theta\|_{L^\infty_T H^2_x}^2+\|\nabla\theta\|_{L^2_T H^2_x}^2
\le C(\kappa_\ast,\|a\|_{C^2}, \kappa^\ast) \|\theta_0\|^2_{H^2}(1+\|\nabla\theta_0\|^2_{L^2})^2
\\
&\qquad\qquad\qquad\times
 \exp\Bigl(C(\kappa_\ast, \|a\|_{C^1})(\|\nabla u\|_{L^2_T L^2_x}^2+\|u\|_{L^4_T L^4_x}^4+\|\nabla\theta \|_{L^4_T L^4_x}^4)\Bigr) .
\end{align*}
\end{paragraph}

    \subsubsection{\texorpdfstring{Case  $(\theta_0,u_0)\in H^{1}(\R^2)\times (H^s(\R^2))^2$, $s\in (0,2]$  }{}}\label{subs:vertical}
      In this subsection we are going to sketch the proof of the $H^s$, $s\in (0,2)$-Estimate \eqref{u:Hsnu} for the divergence-free vector field $u$ of the unique solution $(\theta, u)$ to the Boussinesq equations \eqref{BS}, under the assumption that $\theta_0\in H^1(\R^2)$, following the procedure in Subsection \ref{subs:parabolic}.
      
We deal with the  endpoint case  $(\theta_0, u_0)\in H^1(\R^2)\times (H^2(\R^2))^2$ similarly as for the endpoint case above.
      
      \begin{paragraph}{Case $(\theta_0,u_0)\in H^{1}(\R^2)\times (H^s(\R^2))^2$, $s\in (0,2)$}
       Recall \eqref{LerayProjector} for the definition of the Leray-Helmholtz projector $\PP$ such that
       $$
       \PP u=u,\quad \PP \nabla\Pi=0.
       $$
       We apply $\PP$ to the velocity equation $\eqref{BS}_2$ to arrive at
       \begin{equation}\label{Pu}
       \d_t u+\PP(u\cdot\nabla u)-\PP\di(\mu Su)=\PP(\theta \vec{e_2}).
       \end{equation}
       We  apply $\Delta_j$ to the above equation \eqref{Pu} to arrive at the equation for $u_j:=\Delta_j u$
        \begin{equation}\label{equation:uj}
        \d_tu_j+\PP u\cdot\nabla u_j-\PP\di(\mu Su_j)
        =\PP [u,\Delta_j]\cdot\nabla u-\PP \div([\mu,\Delta_j] Su) +\PP(\theta_j \vec{e_2}).    
        \end{equation}
  We take the $L^2(\R^2)$-inner product between \eqref{equation:uj} and the divergence-free dyadic piece $u_j=\PP u_j$ and follow  the similar argument as to arrive at  \eqref{theta:j:l2:3}, to deduce
    \begin{equation}\label{u:s}
       \frac{d}{dt}\|u_j\|_{L^2_x}+2^{2j}\|u_j\|_{L^2_x}\le C(\mu_\ast) \left(\|[u,\Delta_j]\cdot\nabla u\|_{L^2_x}+2^j\|[\mu,\Delta_j]\nabla u\|_{L^2_x}+\|\theta_j\|_{L^2_x}\right).
    \end{equation}
    By use of the commutator estimate \eqref{commutator:1} in Lemma \ref{Liao} again, we have the following commutator estimates as in \eqref{comm:est}:
     \begin{equation*}\begin{split}
     &\|[u,\Delta_j]\nabla u\|_{L^{2}_x} \le C l_j 2^{j(1-s)}\|\nabla u\|_{L^2_x}\|\nabla u\|_{H^{s-1}_x},
     \hbox{ for } s\in (0,2),\\ 
 &2^j\|[\mu,\Delta_j]\nabla u\|_{L^2_x}
     \le C l_j 2^{j(1-s)}\|\nabla \mu\|_{H^{\nu}_x}\|\nabla u\|_{H^{s-\nu}_x},   
     \hbox{ for } \nu\in (-1,1),\, s\in (-1,\nu+1).
     \end{split}\end{equation*}   
     By virtue of the composition estimate \eqref{composition} in  Lemma \ref{Liao}: 
     $$\|\nabla\mu\|_{H^\nu_x}\leq C(\|b\|_{C^{[\nu]+2}}, \|\theta\|_{H^1})\|\nabla\theta\|_{H^\nu_x},$$
      we derive similar as \eqref{psi:j:s:t} that, for $0<s<\nu+1<2$,
\begin{equation}\label{u:j:s:t}
\begin{aligned}
   & 2^{2js}\|u_j\|_{L^\infty_T L^2_x}^2+2^{2j(s+1)}\|u_j\|_{L^2_T L^2_x}^2 \le 2^{2js}\|(u_0)_j\|_{L^2_x}^2
   \\
& 
   +C(\mu_\ast)\int^T_0 2^{2j(s-1)}\|\theta_j\|_{L^2_x}^2 dt 
   +C(\mu_\ast,s,\nu,\|b\|_{C^{[\nu]+2}}, \|\theta\|_{L^\infty_T H^1_x}) (l_j)^2\times
   \\
   &\qquad
   \times \int^T_0 \Bigl( \|\nabla u\|_{L^2_x}^2\|\nabla u\|_{H^{s-1}_x}^2+\|\nabla \theta\|_{H^{\nu}_x}^2\|\nabla u\|_{H^{s-\nu}_x}^2 \Bigr)\dt.
\end{aligned}    
\end{equation}
By  the interpolation inequality \eqref{interpolation}:
$$
\|\nabla \theta\|_{H^{\nu}_x}\|\nabla u\|_{H^{s-\nu}_x}
\leq C\|\nabla\theta\|_{L^2_x}^{1-\nu}\|\nabla\theta\|_{H^1_x}^\nu \|\nabla u\|_{H^{s-1}_x}^{\nu}\|\nabla u\|_{H^{s}_x}^{1-\nu},
\quad \hbox{ for }\nu\in(0,1),
$$ 
and the  $L^2$-Estimate \eqref{energy:u:1}: 
\begin{equation}\label{u:L2}
\|u\|_{L^\infty_T L^2_x}^2+\|\nabla u\|_{L^2_T L^2_x}^2\leq C(\mu_\ast) (\|u_0\|_{L^2_x}^2+T\|\theta_0\|_{L^2_x}^2),
\end{equation} we arrive at the following by Young's inequality
\begin{align*}
&\|u\|_{L^\infty_T H^s_x}^2+\|\nabla u\|_{L^2_T H^{s}_x}^2
 \leq C(\mu_\ast)\Bigl( \|u_0\|_{H^s_x}^2+T\|\theta_0\|_{  L^2_x}^2+\|\theta\|_{L^2_T H^{s-1}_x}^2\Bigr) 
\\
&\qquad +C(\mu_\ast,s,\nu,\|b\|_{C^2}, \|\theta\|_{L^\infty_T H^1_x}) \int^T_0 (\|\nabla u\|_{L^2_x}^2
+ \|\nabla\theta\|_{H^1_x}^2)\|\nabla u\|_{H^{s-1}_x}^2 dt,
\end{align*}
which, together with Gronwall's inequality, implies \eqref{u:Hsnu}.
\end{paragraph}

\bigbreak

\begin{paragraph}{Endpoint case $(\theta_0, u_0)\in H^1(\R^2)\times (H^2(\R^2))^2$}
We recall the $u$-equation \eqref{Pu} where
$$
\div(\mu Su)=\mu\Delta u+\nabla \mu\cdot Su.
$$
We test \eqref{Pu} by  the divergence-free vector field $\Delta^2 u$, to arrive at
\begin{align*}
&\frac12\frac{d}{dt}\int_{\R^2}|\Delta u|^2\dx
+\int_{\R^2}\mu|\nabla\Delta u|^2\dx
\\
&= \int_{\R^2} \Bigl( -u\cdot\nabla u \Delta^2 u+\nabla\mu\cdot Su\cdot\Delta^2 u
-\nabla\mu\cdot\nabla\Delta u\cdot\Delta u +\Delta\theta\Delta u_2 \Bigr)\dx.
\end{align*}
By use of the embedding  $L^1(\R^2)\hookrightarrow H^{-1}(\R^2)$ again, the righthand side can be bounded by
\begin{align*}
&C(\|\nabla u\|_{L^4_x}^2+\|u\|_{L^4_x}\|\nabla^2u\|_{L^4_x}+\|\nabla^2\mu\|_{L^2_x} \|\nabla u\|_{H^1_x}+\|\nabla\mu\|_{L^4_x} \|\nabla^2u\|_{L^4_x}) \|\nabla\Delta u\|_{L^2_x} 
\\
&+\|\Delta\theta\|_{L^2_x}\|\Delta u\|_{L^2_x}.
\end{align*}
Thus we have the following a priori $\dot H^2_x$-Estimate for $u$ and any positive time $T>0$ by Young's inequality and Gronwall's inequality 
\begin{align*}
 &\|\Delta u(T)\|_{L^2_x}^2+\| \nabla\Delta u\|_{L^2_T L^2_x}^2
 \\
&\leq  C(\mu_\ast)(\|\Delta u_0\|_{L^2_x}^2+\|\nabla u\|_{L^4_T L^4_x}^4+\|\nabla^2\mu\|_{L^2_T L^2_x}^2\|\nabla u\|_{L^\infty_T L^2_x}^2+\|\Delta\theta\|_{L^2_T L^2_x}\|\Delta u\|_{L^2_T L^2_x} )
\\
&\qquad\times \exp\Bigl(C(\mu_\ast)(\|(u, \nabla \mu)\|_{L^4_T L^4_x}^4+\|\nabla^2\mu \|_{L^2_T L^2_x}^2)\Bigr),
\end{align*} 
which gives \eqref{u:Hs2}.
\end{paragraph}


\subsubsection{\texorpdfstring{Case $(\theta_0, u_0)\in H^{s}(\R^2)\times (H^{s-2}(\R^2))^2$ 
or $ H^{s-1}\times (H^s(\R^2))^2$, $s>2$  }{}}\label{subs:s:2}
We are going to use the estimates in the high regularity regime in Lemma \ref{Liao} to derive the $H^s$-estimates \eqref{theta:Hs23}-\eqref{theta:Hs3}-\eqref{u:Hs23}-\eqref{u:Hs3} in Remark \ref{rem:Hs}.
Let $(l'_j)_{j\ge 0}$ be a normalised sequence in $\ell^2(\N)$ such that $l'_j\geq 0$ and $\sum_{j\ge 0}(l'_j)^2=1$.  

\begin{paragraph}{Case  $(\theta_0,u_0)\in H^{s}(\R^2)\times (H^{s-2}(\R^2))^2$, $s>2$}
We can view the transport term $u\cdot\nabla\theta$ simply as a source term of the $\theta$-equation:
$$
\d_t\theta-\div(\kappa\nabla\theta)=-u\cdot\nabla\theta
$$
Then the preliminary estimate for $\theta_j=\Delta_j\theta$ in \eqref{theta:j} can be replaced by  
    \begin{equation*}\label{theta:j+}
       \frac{d}{dt}\|\theta_j\|_{L^2_x}+2^{2j}\|\theta_j\|_{L^2_x}
       \le C(\kappa_\ast) \left(\|(u\cdot\nabla\theta)_j\|_{L^2_x}+2^j\|[\kappa,\Delta_j]\nabla\theta\|_{L^2_x}\right),\quad j\geq 0.
    \end{equation*} 
We apply Lemma \ref{Liao} to derive the following  estimates for $s>1$:
\begin{align*} 
& \|(u\cdot\nabla\theta)_j\|_{L^2_x}
     \le C l'_j 2^{j(1-s)} (\|  u\|_{L^\infty_x}\|\nabla \theta\|_{H^{s-1}_x}+\|  u\|_{H^{s-1}_x}\|\nabla \theta\|_{L^\infty_x}), \\
     &  2^j\|[\kappa,\Delta_j]\nabla\theta\|_{L^2_x}
    \le C l'_j 2^{ j(1-s)} ( \|\nabla \kappa\|_{L^\infty_x}\|\nabla \theta\|_{H^{s-1}_x}+\|\nabla\kappa\|_{H^{s-1}_x}\|\nabla\theta\|_{L^\infty_x}). 
\end{align*}  
        Therefore we have the following estimate similarly as in \eqref{theta:j:s:t} for $s\in (2,3)$ by virtue of $H^{s-1}(\R^2)\hookrightarrow L^\infty(\R^2)$:
\begin{equation*}
\begin{aligned}
   & 2^{2js}\|\theta_j\|_{L^\infty_T L^2_x}^2+2^{2j(s+1)}\|\theta_j\|_{L^2_T L^2_x}^2 \le 2^{2js}\|(\theta_0)_j\|_{L^2_x}^2 \\
    &\quad+ C(\kappa_\ast, s) (l'_j)^2\int^T_0 \| u\|_{H^{s-1}_x}^2\|\nabla\theta\|_{H^{s-1}_x}^2\dt
    \\
    &\quad
    + C(\kappa_\ast,  a,  \|\theta\|_{L^\infty_T L^\infty_x})(l_j')^2 \int^T_0 \|\nabla \theta\|_{L^\infty_x}^2\|\nabla\theta\|_{H^{s-1}_x}^2 dt,\quad j\geq 0,
\end{aligned}    
\end{equation*}
which, together with    the $L^2$-estimate \eqref{theta:L2} and the Gronwall's inequality, implies \eqref{theta:Hs23}.

For $s\geq 3$, we make use of the following commutator estimate 
\begin{align}\label{comm:est:3}  
      \|[u,\Delta_j]\nabla\theta\|_{L^{2}_x}
     &\le C l'_j 2^{j(1-s)}( \|\nabla u\|_{L^\infty_x}\|\nabla \theta\|_{H^{s-2}_x}+\|\nabla u\|_{H^{s-2}_x}\|\nabla \theta\|_{L^\infty_x}),
     \end{align}  
     such that the estimate \eqref{theta:Hs3} follows.
     \end{paragraph}

\begin{paragraph}{Case  $(\theta_0,u_0)\in H^{s-1}(\R^2)\times (H^{s}(\R^2))^2$, $s>2$}
We recall the preliminary estimate for $u_j$ in \eqref{u:s}. 
 We apply Lemma \ref{Liao} to derive  the following  commutator estimates  for $s\in (2,3)$ and $\nu\in (s-2,1)\subset(0,1)$
    \begin{equation*}\label{comm:est:23}\begin{split}
    & \|[u,\Delta_j]\nabla u\|_{L^{2}_x}
     \le C l_j 2^{j(1-s)}\|\nabla u\|_{H^{\nu}_x}\|\nabla u\|_{H^{s-1-\nu}_x}, \\
     &  2^j\|[\mu,\Delta_j]\nabla u\|_{L^2_x}
    \le C l'_j 2^{ j(1-s)} ( \|\nabla \mu\|_{L^\infty_x}\|\nabla u\|_{H^{s-1}_x}+\|\nabla\mu\|_{H^{s-1}_x}\|\nabla u\|_{L^\infty_x}),
        \end{split}\end{equation*}
        which implies then
 \begin{align*} 
   & 2^{2js}\|u_j\|_{L^\infty_T L^2_x}^2+2^{2j(s+1)}\|u_j\|_{L^2_T L^2_x}^2 \le 2^{2js}\|(u_0)_j\|_{L^2_x}^2
   +C(\mu_\ast)\int^T_0 2^{2j(s-1)}\|\theta_j\|_{L^2_x}^2 dt \\
    & + C(\mu_\ast,s,\nu) (l_j)^2 \int^T_0\|\nabla u\|_{H^{\nu}_x}^2\|\nabla u\|_{H^{s-1-\nu}_x}^2\dt
    \\
&    +C(\mu_\ast,s,\|b\|_{C^{[s]+1}}, \|\theta\|_{L^\infty_TH^1_x}) (l'_j)^2\int^T_0  \|\nabla \theta\|_{L^\infty_x}^2\|\nabla u\|_{H^{s-1}_x}^2 
    +\|\nabla \theta\|_{H^{s-1}_x}^2\|\nabla u\|_{L^\infty_x}^2 dt,
    \quad j\geq 0.
 \end{align*} 
This, together with the $L^2$-estimate \eqref{u:L2} and Sobolev's embedding $H^{s-1}(\R^2)\hookrightarrow L^\infty(\R^2)$, implies \eqref{u:Hs23} where $\nu\in (0,1)$ is taken to be a small  constant bigger than $s-2$.

For $s\geq 3$, we use the commutator estimate \eqref{comm:est:3} with $\theta$ replaced by $u$, to arrive at \eqref{u:Hs3}.
\end{paragraph}

\end{document}